\documentclass[11pt,fullpage]{article}

\usepackage{graphicx}
\usepackage{latexsym}
\usepackage{color}
\usepackage{multirow}

\usepackage{amsmath}
\usepackage{graphicx}
\usepackage{amssymb}
\usepackage{amsthm}
\usepackage{amsfonts}

\newcommand {\bb}[1]{{\mathbb #1}}

\renewcommand{\baselinestretch} {1.3}

	\makeatletter 
\def\singlespace{\def\baselinestretch{1}\@normalsize}

\newtheorem{theorem}{Theorem}
\newtheorem{lemma}{Lemma}

\newtheorem{corollary}{Corollary}

\newtheorem{remark}{Remark}

\newcommand{\be}{\begin{equation}}
\newcommand{\ee}{\end{equation}}
\newcommand{\beqn}{\begin{eqnarray}}
\newcommand{\eeqn}{\end{eqnarray}}
\newcommand{\beqns}{\begin{eqnarray*}}
\newcommand{\eeqns}{\end{eqnarray*}}

\newcommand{\fr}[1]{(\ref{#1})} 
\newcommand{\lkr}{\left(}  
\newcommand{\rkr}{\right)} 
\newcommand{\lkv}{\left[}  
\newcommand{\rkv}{\right]} 
\newcommand{\lfi}{\left\{} 

\newcommand{\EE}{\ensuremath{{\mathbb E}}}

\newcommand{\ZZ}{\ensuremath{{\mathbb Z}}}

\newcommand{\wind}{\mbox{wind}}
\newcommand{\Tr}{\mbox{Tr}}

\newcommand{\pen}{\mbox{pen}}
\newcommand{\Var}{\mbox{Var}}
\newcommand{\AIF}{\mbox{AIF}}

\newcommand{\om}{\omega}

\newcommand{\bA}{\mbox{\boldmath $A$}}
\newcommand{\bU}{\mbox{\boldmath $U$}}
\newcommand{\bI}{\mbox{\boldmath $I$}}
\newcommand{\bG}{\mbox{\boldmath $G$}}
\newcommand{\bQ}{\mbox{\boldmath $Q$}}
\newcommand{\bD}{\mbox{\boldmath $D$}}

\newcommand{\bz}{\mbox{\boldmath $z$}}
\newcommand{\bs}{\mbox{\boldmath $s$}}
\newcommand{\bt}{\mbox{\boldmath $t$}}
\newcommand{\bq}{\mbox{\boldmath $q$}}
\newcommand{\bg}{\mbox{\boldmath $g$}}
\newcommand{\bof}{\mbox{\boldmath $f$}}
\newcommand{\by}{\mbox{\boldmath $y$}}
\newcommand{\bh}{\mbox{\boldmath $h$}}
\newcommand{\bzero}{\mbox{\boldmath $0$}}
\newcommand{\bxi}{\mbox{\mathversion{bold}$\xi$}}
\newcommand{\boeta}{\mbox{\mathversion{bold}$\eta$}}
\newcommand{\bzeta}{\mbox{\mathversion{bold}$\zeta$}}
\newcommand{\bepsilon}{\mbox{\mathversion{bold}$\epsilon$}}

\newcommand{\bPhi}{\mbox{\mathversion{bold}$\Phi$}}
\newcommand{\bOm}{\mbox{\mathversion{bold}$\Omega$}}

\newcommand{\vf}{ \vec{\bof}}
\newcommand{\vbof}{ \vec{\bof}}
\newcommand{\vxi}{\vec{\bxi}}
\newcommand{\vz}{\vec{\bz}}
\newcommand{\vq}{\vec{\bq}}
\newcommand{\veta}{\vec{\boeta}}

\newcommand{\whm}{\widehat{m}}

\newcommand{\mathM}{{\mathcal M}}

\graphicspath{%
    {converted_graphics/}
    {/}
}
\begin{document}

\title{\Large{\bf Laplace deconvolution and its application to Dynamic Contrast Enhanced imaging}}

\author{
\large{{\sc Fabienne Comte$^1$}, {\sc Charles-Andr\'{e} Cuenod$^{1,2}$ },   {\sc Marianna Pensky$^3$}}\\  
\large{and {\sc Yves Rozenholc$^1$}}
  \\ \\
   Universit\'{e}  Paris Descartes$^1$, European Hospital George Pompidou$^2$\\ and  University of Central Florida$^3$ }

\date{}

\maketitle

\vspace{1cm}
\begin{abstract}
In the present paper we consider the problem of Laplace deconvolution 
with noisy discrete observations. The study is motivated by 
Dynamic Contrast Enhanced imaging  using a bolus of contrast agent, 
a procedure  which allows considerable improvement in  {evaluating}
the quality of a vascular network and its permeability and is widely used 
 in medical assessment  of brain flows or cancerous tumors.
Although the study is motivated by medical imaging application, 
we obtain a solution of a general problem of Laplace deconvolution  based on noisy data
which appears in many different contexts.
We propose a new method for Laplace deconvolution which is based on 
expansions of the convolution kernel, the unknown function and the observed 
signal over   Laguerre functions basis. 
The expansion results in a small system of linear equations with the matrix of the 
system being triangular and Toeplitz. 
The number  $m$ of the terms in  the expansion  of the estimator is controlled via complexity penalty.
The advantage of this methodology is that it leads to very fast computations,
does not require exact knowledge of the kernel 
and produces no  boundary effects due to extension at zero and cut-off at $T$.
The technique  leads to an estimator with the risk within a logarithmic factor 
of $m$ of the oracle risk under no assumptions on the model and within a 
constant factor of the oracle risk under mild assumptions.
 The methodology is illustrated by a finite sample simulation study 
which includes an  example of the kernel obtained in the real life DCE experiments.
Simulations confirm that the proposed technique is fast, efficient, accurate,  
usable from a practical point of view  and competitive.
\end{abstract}


\vspace{1cm}

\bibliographystyle{plain}

\pagestyle{plain}

\section{Introduction}
\label{sec:intro}
\setcounter{equation}{0}

Cancers and vascular diseases present major public health concerns. Considerable improvement in assessing
the quality of a vascular network and its permeability have been achieved through Dynamic Contrast
Enhanced   (DCE) imaging using a bolus of contrast agent at high frequency such as 
Dynamic Contrast Enhanced Computer Tomography (DCE-CT),
Dynamic Contrast Enhanced Magnetic Resonance Imaging  (DCE-MRI) and Dynamic Contrast Enhanced Ultra Sound (DCE-US). 
Such techniques  are   widely used   in medical assessment  of brain flows or cancerous tumors
(see, e.g., Cao {\it et al.}, 2010; Goh {\it et al.}, 2005; Goh and Padhani, 2007;
Cuenod {\it et al.}, 2006; Cuenod {\it et al.}, 2011;  Miles, 2003;  Padhani and Harvey,  2005 and Bisdas
{\it et al.}, 2007). This imaging procedure has great potential for cancer detection and
characterization, as well as for monitoring \textit{in vivo} the effects of treatments.
It is also used, for example, after a stroke for prognostic purposes or for occular blood flow 
evaluation.

\begin{figure}
\[ \includegraphics[height=8cm]{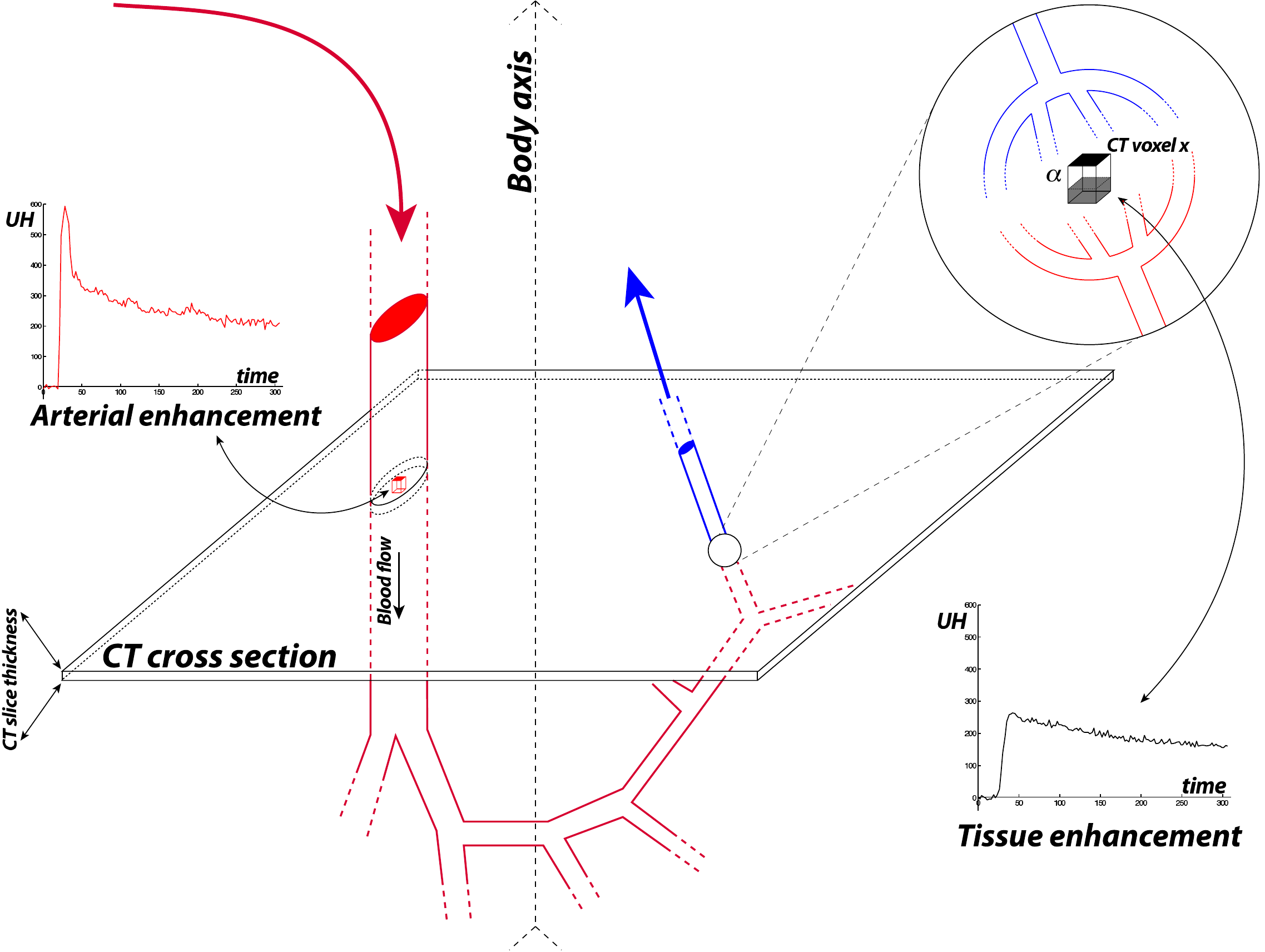}\]
\caption{\label{fig:DCE-CT} {\bf DCE-CT experiment and contrast agent circulation.} 
The patient body is materialized by the mixed arrow. }
\end{figure}

As an example, below we consider a DCE-CT experiment that  follows  the diffusion  of a 
bolus of a contrast agent injected into a vein.  At the microscopic level,
for a given voxel of interest   having unit volume,   the number of arriving particles at time $t$ is given by 
$\beta\, \AIF(t)$, where the Arterial Input Function (AIF) measures  concentration within a unit 
volume voxel inside the aorta  and $\beta$ is a   proportion  of the AIF which enters the tissue voxel.
Denote the number of particles in the voxel   at time $t$ by $Y(t)$ and the random lapse of time 
during which a particle sojourns in the voxel by $S$. Assuming sojourn times for different particles 
to be i.i.d. with   c.d.f. $F$, one obtains the following equation for the average number of 
contrast agent particles at the moment $t$
$$
{\bb E} Y(t) =  \hspace{-3mm}\underbrace{\int_0^t \beta\, \AIF (t-\tau)\,d\tau}_{\text{arrived before time $t$}}
\hspace{-1mm}- \underbrace{\int_0^t  \beta\, \AIF (t-\tau)\,P(S \leq \tau)\,d\tau}_{\text{left before time $t$}} 
= \int_0^t  \beta\, \AIF (t-\tau)(1-F(\tau)) d\tau,
$$
where the expectation is taken under the unknown distribution of the sojourn times.
In reality, one does not know ${\bb E} Y(t) $ and has discrete noisy observations
$$
Y(t_i) = {\bb E} Y(t_i) + \sigma \epsilon_i.
$$

Medical doctors are interested in a reproducible quantification of the blood flow inside 
the tissue which is characterized by $f(t)= \beta(1-F(t))$
since this quantity is independent of the concentration of particles of contrast agent 
within a unit volume voxel inside the aorta  described by $\AIF(t)$.
The sequential imaging acquisition is illustrated by Figure \ref{fig:DCE-CT}. 
The contrast agent arrives with the oxygenated blood through the aorta (red arrow) 
where its concentration, AIF, within unit volume voxel   is first measured when it passes 
through the CT cross section (red box). 
Subsequently, the contrast agent enters   the arterial system, and it is assumed 
that its concentration does not change during this phase. The exchange within the tissue 
of both oxygen and contrast agent occurs after the arterial phase and the concentration 
of contrast agent during this exchange is measured in all tissue voxels (grey voxel in the zoom) 
inside the CT cross section. Later the contrast agent returns to the venous system with the 
de-oxygenated blood (blue arrow).

To complete description of this experiment, one has to take into account that  there is a  
delay $\delta$ between the measurement of the contrast agent concentration inside the aorta  
(first cross of the CT section) and its arrival inside the tissue. This leads to the following complete model:
\be \label{eq:AIF model}
Y(t_i) = \int_0^{t_i-\delta} \beta \AIF(t_i - \tau) (1-F(\tau)) d\tau + \sigma \epsilon_i, \;\;\; i=1,...,n.
\ee
The value of delay $\delta$ can be  measured with a small error using the decay between the jumps after the injection of the contrast agent  inside the aorta and the tissue. Unfortunately, evaluation of the proportion $\beta$ is a much harder task which is realized  with a larger error.
In the spirit of   complete model \fr{eq:AIF model}  for DCE-CT experiments, one can consider a more general model of
Laplace convolution equation based on noisy observations which presents  a necessary theoretical step 
before obtaining medical answers provided by model \fr{eq:AIF model}.
\\

Indeed, for a known value  of  $\delta$, equation \fr{eq:AIF model} reduces   to a noisy version of a Laplace
convolution equation
\be \label{eq:model}
y(t_i) = \int_0^{t_i} g(t_i - \tau) f(\tau) d\tau + \sigma \epsilon_i, \;\;\; i=1, \ldots, n,
\ee
where function $g$ is considered to be known, $f$ is a function of interest,
measurements $y(t_i)$ are taken at points $0 \leq t_1 \leq ... \leq t_n \leq T < \infty$ and 
$\epsilon_i$ are i.i.d. $N(0,1)$. The corresponding noiseless version of this equation can be written as 
\be \label{eq:Volterra}
q(t) = \int_0^t g(t - \tau) f(\tau) d\tau, \;t \geq 0.
\ee

Formally, by setting $g(t) = f(t) \equiv 0$ for $t<0$, equation
\fr{eq:Volterra} can be viewed as a particular case of the
Fredholm convolution equation
\be \label{eq:Fredholm}
q(t) = \int_{a}^b  g(t  - \tau) f(\tau) d\tau,
\ee
where $a= -\infty$ and $b=\infty$ for Fourier convolution on a real line and 
$-\infty < a < b < \infty$ for circular convolution.
Discrete stochastic version of equation \fr{eq:Fredholm}
\be \label{eq:fourier_deconv}
y(t_i) = \int_{a}^b g(t_i - \tau) f(\tau) d\tau + \sigma \epsilon_i, \;\;\;\; i=1,...,n,
\ee
known also as Fourier deconvolution problem,
has been extensively studied in the last thirty years
(see, for example, Carroll and Hall, 1988; Comte, Rozenholc and Taupin, 2007; 
Delaigle, Hall and Meister, 2008;
Diggle and Hall, 1993; Fan, 1991; Fan and Koo, 2002; Johnstone {\it et al.},
2004; Pensky and Vidakovic, 1999;  Stefanski and Carroll, 1990, among others).
However, such an approach to solving (\ref{eq:Volterra}) and (\ref{eq:model})
is very misleading.

Indeed, since one does not have data outside the interval $[0,T]$ and since function 
$f(t)$ may not vanish fast enough as $t \to \infty$, one cannot apply Fourier transform 
on the whole real line since Fourier transform is defined for only integrable or square integrable functions.
Application of the discrete Fourier transform (DFT) on the finite interval $[0,T]$ is useless since
the kernel $g$ is not periodic. Consequently, convolution in equation \fr{eq:Volterra}
is not circular and, hence, it is not converted into a product by DFT.

The  issue of having measurements only on the part $t \leq T$ of half line $(0,\infty)$
does not affect the Laplace deconvolution since it exhibits {\bf causality} property:
the values of $q(t)$ for $0\leq t \leq T$ depend on values of $f(t)$ for $0\leq t \leq T$
only and vice versa.

The mathematical theory of (noiseless) convolution type Volterra equations
is well developed (see, e.g., Gripenberg  {\it et al.} 1990)
and the exact solution of \fr{eq:Volterra} can be obtained through
Laplace transform.  However, direct application of Laplace transform for
discrete measurements faces serious conceptual and numerical problems.
The inverse Laplace transform is usually  found by application of tables of inverse Laplace
transforms, partial fraction decomposition or series expansion 
(see, e.g., Polyanin and Manzhirov, 1998), neither of which is
applicable in the case of the discrete noisy version of Laplace deconvolution.
Only few applied mathematicians and researchers in natural sciences took an effort
to solve the problem using discrete measurements in the left hand side of \fr{eq:Fredholm}.
Since the problem arises in medical imaging, few scientists put an effort to solve equation \fr{eq:AIF model} using 
singular value decomposition (SVD) with the subsequent application of Tikhonov regularization
(see, e.g., Axel (1980), Ostergaard {\it et al.} (1996) and   an extensive review in
Fieselmann {\it et al.} (2011)). In fact, SVD has been widely used in the context of DCE imaging since mid-nineties. 
The technique, however, is very computationally unstable, especially, in the presence of recirculation 
of contrast agent. For this reason, SVD has been  mostly used in the simplified framework  of brain imaging due to the presence of white barrier
which prevents circulation of contrast agent outside blood vessels. 
Ameloot and Hendrickx  (1983) applied Laplace deconvolution for the
analysis of  fluorescence curves and used a parametric presentation of the solution $f$ as a sum of exponential functions
with parameters evaluated  by minimizing discrepancy with the right-hand side.
In a somewhat similar manner,  Maleknejad {\it et al.} (2007) proposed to expand
the unknown solution over a wavelet basis and find the coefficients via the least squares algorithm.  
Lien {\it et al.} (2008), following Weeks  (1966), studied numerical inversion of the Laplace transform using Laguerre functions.
Finally, Lamm (1996) and Cinzori and Lamm (2000) used discretization of
the equation \fr{eq:Volterra} and  applied various versions of the Tikhonov regularization  technique.
However, in all of the above papers, the noise in the measurements  was either ignored
or treated as deterministic. The presence of random noise in (\ref{eq:model})
makes the problem even more challenging.

For   the  reasons listed above, estimation of $f$ from discrete noisy
observations $y$ in \fr{eq:model} requires extensive investigation.
Unlike Fourier deconvolution that has been intensively studied in statistical
literature (see references above), Laplace deconvolution
received very little  attention within statistical framework.
To the best of our knowledge, the only paper which tackles the problem is
Dey, Martin and Ruymgaart (1998) which considers a noisy version of Laplace deconvolution
with a very specific kernel of the form $g(t)=be^{-at}$. The authors use the fact that, in this case,
the solution of the equation \fr{eq:Volterra}  satisfies a particular linear differential equation and,
hence, can be recovered using   $q(t)$ and its derivative $q' (t)$.
For this particular kind of kernel, the authors   derived   convergence
rates for the  quadratic risk of the proposed estimators, as $n$ increases,
under the assumption that  the $s$-th derivative of
$f$ is   continuous on $(0,\infty)$. However, in Dey, Martin and Ruymgaart (1998)
it is assumed that data are available on the whole positive half-line (i.e.  $T= \infty$)
and that $s$ is known (i.e., the estimator is not adaptive).

Recently, Abramovich {\it et al.} (2012) studied the problem of Laplace deconvolution 
based on discrete noisy data. The idea of the method is to reduce the problem to 
estimation of the unknown regression function, its derivatives and, possibly, 
some linear functionals of these derivatives. The estimation is carried out using kernel method with the 
 Lepskii technique for the choice of the bandwidth (although it is mentioned in the paper that other 
methodologies for the choice of bandwidth can also be applied). The method has an advantage of reducing 
a new statistical problem to a well studied one. However, the shortcoming of the technique is that 
it requires meticulous  boundary correction and is strongly dependent on the knowledge of the kernel $g$.
Indeed, small change in the kernel may produce significant changes in the expression for the estimator.

In the present paper we suggest a method which is designed to overcome limitations of 
the previously developed techniques. The new methodology is based on expansions of the kernel, 
unknown function $f$ and the right-hand side in equation   \fr{eq:model} over the   Laguerre functions basis. 
The expansion results in a small system of linear equations with the matrix of the 
system being triangular and Toeplitz. 
The number of the terms in  the expansion  of the estimator is controlled via complexity penalty.
The advantage of this methodology is that it leads to very fast computations
and produces no  boundary effects due to extension at zero and cut-off at $T$.
The technique does not require exact knowledge of the kernel since it is represented
by its Laguerre coefficients only and leads to an estimator with the risk within a logarithmic factor 
of $m$ of the oracle risk under no assumptions on the model and within a constant factor of the oracle risk
under mild assumptions. Another merit  of the new methodology includes the fact that, since the unknown 
functions are represented by a small number of Laguerre coefficients, it is easy to 
cluster or classify them for various groups of patients.
Simulation study shows that the method is very accurate and stable and easily outperforms 
SVD and kernel-based technique of Abramovich,   Pensky, and Rozenholc  (2012).

The rest of the paper is organized as follows. In Section~\ref{sec:Laguerre_functions} 
we derive the system of equations resulting from expansion of the functions over the Laguerre
basis, study the effect of discrete, possible irregularly spaced data and introduce 
selection of model size via penalization. Corollary~\ref{cor:risk_bound}
indeed confirms that the risk of the penalized estimator lies within a logarithmic factor of $m$ of the
minimal risk. In Section~\ref{sec:asymp_risk} we obtain asymptotic upper bounds for the risk and 
prove  the risk lies within a constant factor of an oracle risk. The proof of this fact rests on
nontrivial facts of the theory of Toeplitz matrices. Section~\ref{sec:simulation}
provides a finite sample simulation studies. Finally, Section~\ref{sec:discussion}
discusses results obtained in the paper. Section~\ref{sec:proofs}
contains proofs of the results in the earlier sections.

\section{Laplace deconvolution via expansion over Laguerre functions basis }
\label{sec:Laguerre_functions}
\setcounter{equation}{0}

\subsection{Relations between coefficients of the Laguerre expansion}
\label{sec:linear_system}

One of the possible solution of the problem \fr{eq:model} is to use Galerkin method with
the basis represented by a system of Laguerre functions. Laguerre functions are defined as
\be \label{eq:Laguerre_func}
\phi_k(t) = \sqrt{2a} e^{-at} L_k(2at),\ \ k=0,1, \ldots,
\ee
where $L_k(t)$ are Laguerre polynomials (see, e.g., Gradshtein and Ryzhik (1980))
$$
L_k(t) = \sum_{j=0}^k (-1)^j {k \choose j} \frac{t^j}{j!},\ \ \ t \geq 0.
$$
It is known that functions  $\phi_k(\cdot)$, $k=0,1, \ldots$, form an orthonormal basis of the $L^2 (0, \infty)$ space
and, therefore, functions $f(\cdot)$, $g(\cdot)$, $q(\cdot)$ and $y(\cdot)$ can be expanded over this basis with coefficients
$f^{(k)}$, $g^{(k)}$, $q^{(k)}$ and $y^{(k)}$, $k = 0, \ldots, \infty$, respectively.
By plugging these expansions into formula \fr{eq:Volterra}, we obtain the following
equation
\be \label{eq:laguerre1}
\sum_{k=0}^\infty q^{(k)} \phi_k(t) = \sum_{k=0}^\infty \sum_{j=0}^\infty f^{(k)} g^{(j)} \int_0^t \phi_k(x)  \phi_j (t-x) dx.
\ee

It turns out that coefficients of interest $f^{(k)}$, $k=0,1, \ldots,$  can be represented as a solution of
an infinite triangular system of linear equations. Indeed, it is easy to check that (see, e.g., 7.411.4 in Gradshtein and Ryzhik (1980))
$$
\int_0^t \phi_k(x)  \phi_j(t-x)  dx = 2a e^{-at} \int_0^t L_k(2at) L_j(2a(t-x)) dx =
(2a)^{-1/2} \,  [\phi_{k+j}(t) - \phi_{k+j+1} (t)].
$$
Hence,   equation \fr{eq:laguerre1} can be re-written as
$$
\sum_{k=0}^\infty q^{(k)} \phi_k(t) = \sum_{k=0}^\infty \phi_k(t) [(2a)^{-1/2} \, f^{(k)} g^{(0)}  +
\sum_{l=0}^{k-1} (2a)^{-1/2} \, (g^{(k-1)}  - g^{(k-l-1)}) f^{(l)}].
$$
Equating coefficients for each basis function, we obtain an infinite triangular system of linear equations.
In order to use this system for estimating $f$, we define
\be \label{eq:fm}
f_m (x) = \sum_{k=0}^{m-1}  {f}^{(k)} \phi_k (x),
\ee
approximation of $f$ based on the first $m$ Laguerre functions. The following Lemma states
how the coefficients in \fr{eq:fm} can be recovered.

\begin{lemma} \label{lem:triangular_system}
Let $\bof_m$, $\bg_m$ and $\bq_m$ be $m$-dimensional vectors with elements $f^{(k)}$, $g^{(k)}$ and $q^{(k)}$, $k=0,1, \ldots, m-1$,
respectively. Then, for any $m$, one has $\bq_m = \bG_m  \bof_m$ where
$\bG_m$ is the lower triangular Toeplitz matrix  with elements
\be  \label{inf_Toeplitz}
G^{(i,j)}  = \left\{
\begin{array}{ll}
(2a)^{-1/2} \, g^{(0)}, & \mbox{if}\ \ i=j,\\
(2a)^{-1/2} \,  (g^{(i-j)}  - g^{(i-j-1)}), & \mbox{if}\ \ j<i,\\
0, & \mbox{if}\ \ j>i.
\end{array} \right.
\ee
\end{lemma}

Hence,  $f(x)$ can be estimated by
\be \label{eq:f_estimator}
\hat{f}_m (x) = \sum_{k=0}^{m-1} \hat{f}^{(k)} \phi_k (x)
\ee
where  $\hat{\bof}_m = \bG_m ^{-1} \hat{\bq}_m$ and   $\hat{\bq}_m$ is an
unbiased estimator of the unknown vector of coefficients $\bq_m$.

\subsection{Recovering Laguerre coefficients from discrete noisy data}
\label{sec:discrete_data}

Unfortunately, unlike some other linear ill-posed problems, data does not come in the form of
unbiased estimators of the unknown coefficients $q^{(k)}$, $k=0,1, \ldots, m-1$. Below, we examine how the length
of the observation interval  $T$ and spacing of observations $t_i$ in equation \fr{eq:model} affect the system
of equations  in Lemma \ref{lem:triangular_system}.

Let $P:[0,T] \rightarrow [0,T]$ be a function generating observations in \fr{eq:model} such that
$P$ is a continuously differentiable strictly increasing function
\be \label{eq:generate}
P(0)=0,\ \ P(T) =T,\ \ P(t_i) = iT/n,\ \ i=1, \ldots, n.
\ee
Under conditions \fr{eq:generate}, $P$ is a one-to-one function and, therefore, has an inverse $P^{-1}$.

Choose $M$  large enough that the bias in representation \fr{eq:fm}  of $f$ by $f_M$ is very small and
form an $(n\times M)$  matrix $\bPhi_M$ with elements $\Phi^{(i,k)} = \phi_k (t_i)$, $i=1, \ldots,n,$ $k=0, \ldots, M-1$.
Let $\bz_M$ be the $M$-dimensional vector with elements $z^{(j)} = \langle y, \phi_j \rangle$, $j=0, \ldots, M-1$.
Then, it follows  that
$$
y(t_i) = \sum_{l=0}^{M-1} z^{(l)} \phi_l (t_i) = (\bPhi_M \bz_M)^{(i)}, \ \ i=1, \ldots, n.
$$
If $\by$ and $\bh$ are $n$-dimensional vectors with components $y(t_i)$ and $q(t_i)$, $i=1, \ldots,n$,
respectively, then the vectors $\bq_M$  and  $\bz_M$ of the true and the estimated Laguerre coefficients
of $q(x)$ can be represented, respectively,  as
\be \label{eq:rhs_coeffs}
\bq_M =  (\bPhi_M^T \bPhi_M)^{-1} \bPhi_M^T \bh,\ \ \
\bz_M = (\bPhi_M^T \bPhi_M)^{-1} \bPhi_M^T \by.
\ee
Let us examine  matrix $\bPhi_M^T \bPhi_M$. Note that, for any $k$ and $l$,
\beqns
(\bPhi_M^T \bPhi_M)^{(k,l)}  & = & \sum_{i=1}^n \phi_k (P^{-1} (iT/n)) \phi_l (P^{-1} (iT/n))
\approx n T^{-1}\, \int_0^T \phi_k (P^{-1} (x)) \phi_l (P^{-1} (x)) dx \\
& = & n T^{-1}\, \int_0^T \phi_k (t) \phi_l (t) p(t) dt
\eeqns
where $p(t)= P' (t)$. It follows from the above that matrix $\bPhi_M^T \bPhi_M$ should be normalized
by a factor $n^{-1} T$. Indeed, if points $t_i$ are equispaced on the interval $(0,T]$, then,
for $n$ and $T$  large enough, $(\bPhi_M^T \bPhi_M) \approx n T^{-1}\, \bI_M$
where $\bI_M$ is the $M$-dimensional identity matrix. Hence, in what follows, we are going to operate with  matrix
$\bA_M = T  n^{-1} (\bPhi_M^T \bPhi_M)$ and its inverse
\be \label{eq:bOm}
\bOm_M = (\bA_M)^{-1} = n T^{-1} (\bPhi_M^T \bPhi_M)^{-1}.
\ee

Let $\bepsilon$ be the vector with components $\epsilon(t_i)$, $i=1, \ldots, n$, and
$\bxi_M = \sqrt{n/T} \ (\bPhi_M^T \bPhi_M)^{-1} \bPhi_M^T \bepsilon$. Then, the vector $\bof_M$
of the {\bf true} Laguerre coefficients of the unknown function $f$ satisfies the following equation
\be \label{eq:coef_model}
\bz_M = \bG_M \bof_M + \sigma \sqrt{T/n}\ \bxi_M,\ \ \bxi_M \sim N(\bzero, \bOm_M).
\ee
If points $t_i$ are equispaced on the interval $(0,T]$ and both $n$ and $T$ are  large,
then,   in \fr{eq:coef_model},  $\bOm_M \approx \bI_M$.

\subsection{Model selection and oracle risk}
\label{sec:model_selec}

Equation \fr{eq:coef_model} implies that one can estimate unknown vector $\bof_M$ by
$\hat{\bof}_M = \bG_M^{-1} \bz_M$. However, since the value of $M$ is large, the variance of this estimator, 
$$
\EE \|\hat{\bof}_M - \bG_M^{-1} \bq_M \|^2 = \frac{\sigma^2  T}{n}\, \Tr(\bG_M^{-1} \bOm_M \bG_M^{-T}),
$$
where $\bA^{-T} = (\bA^T)^{-1} = (\bA^{-1})^T$ for any matrix $\bA$,  will be too large while the bias in the representation \fr{eq:fm} of $f$ by $f_M$ will be very small.
Hence, in order to balance the bias and the variance components of the error, one needs to
choose the best possible number $m$ of Laguerre functions in the representation \fr{eq:fm} of $f$,
i.e., choose the model size.

In order to achieve 
a required balance between
the bias and the variance components of the error, consider a collection of integer indices
$ {\mathcal M}_n = \{ 1, \ldots, M \}$  
where $M<n$ may depend on $n$ and, for $m \in \mathM_n$, the associated subspaces
$S_m \subseteq R^M$ defined by 
$$
\bt \in S_m\ \ \mbox{if}\ \   \bt = (t^{(0)}, t^{(1)}, \ldots, t^{(m-1)}, 0,0, \ldots, 0)^T.
$$
Let us denote by $\vz_m$, $\vq_m$,  $\vbof_m$, $\vxi_m$ and $\widehat{\vbof}_m $ the $M$-dimensional vectors
where the first $m$ elements coincide with the elements of $m$-dimensional vectors $\bz_m$, $\bq_m$, $\bof_m$,  $\bxi_m$ and $\widehat{\bof}_m$ respectively, and the last $(M-m)$ elements are identical zeros.
For each $m \in \mathM_n$, evaluate 
$$
\widehat{\bof}_m = (\bG_m)^{-1}  \bz_m = \bG_M^{-1} \vz_m 
$$
and denote
\be \label{eq:Qm}
\bQ_m =  \bG_m^{-1} \bOm_m \bG_m^{-T}.
\ee

For the  estimator $\hat{f}_m$ of $f$ given by \fr{eq:f_estimator} with the vector of coefficients
$\widehat{\vbof}_m$,  the bias-variance decomposition of the mean
squared error is of the form
\be \label{eq:bias_var}
\EE \|\hat{f}_m - f\|^2 =   \|f_m - f\|^2 + \sigma^2 T n^{-1}\ \Tr(\bQ_m),
\ee
where the bias term $\|f_m - f\|^2 = \sum_{j = m}^\infty (f^{(j)})^2$ is decreasing and the variance term
$\sigma^2 T n^{-1}\ \Tr(\bQ_m)$ 
 is growing with $m$.  
The smallest possible risk, the so-called {\it oracle } risk,   is obtained by minimizing the right-hand side
of expression \fr{eq:bias_var} with respect to $m$:
\be \label{eq:oracle}
R_{oracle} = \min_m\, \EE \|\hat{f}_m - f\|^2 =
\min_m\, \left[ \|f_m - f\|^2 + \sigma^2 T n^{-1}n \ \Tr(\bQ_m)  \right].
\ee
Hence, the objective is to choose a value of $m$ which delivers an  estimator of the unknown function $f(x)$
with the risk as close to the oracle risk \fr{eq:oracle} as possible.
Since the bias in the right-hand side of expression \fr{eq:bias_var} is unknown, in order to attain this goal, one can use
a penalized version of estimator \fr{eq:f_estimator} as it is described in the next section.

\subsection{Selection of model size via penalization}
\label{sec:model_size}

For any vector $\bt \in R^M$, we define contrast as
\be \label{eq:contrast}
\gamma^2_n (\bt) = \|\bt\|^2  - 2 \langle \bt, \bG_M^{-1} \bz_M \rangle
\ee
and note that for $\bt \in S_m$ one has, thanks to the nul coordinates of $\bt$ and the lower triangular form of $\bG_M$ and $\bG_m$, 
$$\langle \bt, \bG_M^{-1} \bz_M \rangle = \langle \bt, \bG_M^{-1} \vz_m \rangle = \langle \bt, \widehat{\vbof}_m  \rangle.$$
Let $\|\bA\|_2 = \sqrt{\Tr(\bA^T \bA)}$ and $\|\bA\| = \sqrt{\lambda_{\max} (\bA^T \bA)}$ be, respectively,
the Frobenius and the spectral norm  of a matrix $\bA$,
where $\lambda_{\max} (\bU)$ is the largest eigenvalue of $\bU$.
Denote
\be \label{eq:vm_rhom}
v_m^2 = \| \sqrt{\bQ_m} \|^2_2 = \Tr(\bQ_m),\ \ \ \rho_m^2 = \|\sqrt{\bQ_m}\|^2  = \lambda_{\max} (\bQ_m)
\ee
where $\sqrt{\bQ_m}$ is a lower triangular matrix such that $(\sqrt{\bQ_m})^T \sqrt{\bQ_m} = \bQ_m$.
Assume that $\rho_m^2$ grows at most polynomially in $m$, i.e.
there exist positive constants $\alpha$ and $C_\rho$ such that
\be \label{eq:rho_condition}
\rho^2_m \leq C_\rho m^\alpha.
\ee
Choose any constant  $B>0$ and introduce a penalty
\be \label{eq:penalty}
\pen(m) =  4 \sigma^2 T n^{-1} \lkv (1+ B) v_m^2 + (1+ B^{-1}) (2 \alpha + 2) \rho_m^2  \, \log m \rkv.
\ee
For each $m=1,\ldots,M$,  
construct  estimator $\hat{f}_m(x)$ of $f(x)$ of the form \fr{eq:f_estimator}
with coefficients $\widehat{\bof}_m = (\bG_m)^{-1}  \bz_m$, 
the augmented version $\widehat{\vbof}_m$ satisfies 
$$
\widehat{\vbof}_m =\arg\min_{\bt \in S_m}\gamma_n^2(\bt).
$$ 
Now,  choose
$m = \hat{m}$ where
\be \label{eq:hatm}
\hat{m} = \arg \min \left\{m \in {\mathcal M}_n:\  \gamma^2_n (\widehat{\vbof}_m) + {\rm pen}(m) \right\}.
\ee
The  following statement holds.

\begin{theorem} \label{th:risk_bound}
Let condition \fr{eq:rho_condition}  hold for some positive constants $\alpha$ and $C_\rho$.
Then, for any   $B>0$, one has
\be \label{eq:risk_bound}
R(\hat{f}_{\widehat{m}}) := \EE (\|\hat{f}_{\widehat{m}} - f\|^2) \leq  \min_{m \in {\mathcal M}_n}  \left[ 3 \|f_m - f\|^2 +  4 {\rm pen}(m)    +
16 C_\rho^2 \sigma^2 (1 + B^{-1}) \, \frac{T}{mn} \right].
\ee
\end{theorem}

The proof of this and later statements are given in Section \ref{sec:proofs}.\\

Note that the upper bound in Theorem \ref{th:risk_bound} is non-asymptotic and holds
for any values of $T$ and $n$ and any distribution of points $t_i$, $i=1, \ldots, n$.

In order to evaluate   relative precision   of the estimator $\hat{f}_{\widehat{m}}$ constructed above, we shall
compare its risk with the  oracle risk \fr{eq:oracle}.
Since    $\rho_m^2 \leq v_m^2$  for any value of $m$, it follows from Theorem \ref{th:risk_bound}
that,  for any value of $m$, the risk of the estimator $\hat{f}_{\widehat{m}}$ lies within a
logarithmic factor of the oracle risk, i.e., the estimator is optimal within a logarithmic factor of $m$.
In particular, the following corollary   holds.

\begin{corollary} \label{cor:risk_bound}
Under conditions of Theorem \ref{th:risk_bound},
\be \label{eq:risk2}
R(\hat{f}_{\widehat{m}}) \leq 16 [(1+ B) +  (1+ B^{-1}) (2 \alpha + 2)  \log m_0] R_{oracle} +
16 C_\rho^2 \sigma^2 (1 + B^{-1}) \, T  m_0^{-1} n^{-1},
\ee
where $m_0=m_0(n, T)$ is the value of $m$ delivering the minimum in the right-hand side of \fr{eq:risk_bound}.
\end{corollary}

\section{Asymptotic upper bounds for the risk and optimality of the estimator }
\label{sec:asymp_risk}
\setcounter{equation}{0}

\subsection{Assumptions}

Corollary \ref{cor:risk_bound} is valid for any function $g$ and any distribution of sampling points,
hence, it is true in the ``worst case scenario''. In majority of practical situations, however,
$v_m^2$ increases much faster with $m$ than $\rho_m^2$ and the risk of the estimator  $\hat{f}_{\widehat{m}}$
can exceed the oracle risk only by a finite factor independent of $m_0$ and $n$.
In particular, in what follows, we shall show that,
under certain conditions, for $n$   large enough and   $T = T_n$,
the  ratio between  $R(\hat{f}_{\widehat{m}})$ and  $R_{oracle}$
is bounded by a constant independent of $n$.

 For this purpose, assume that  function $g(x)$, its Laplace transform $G(s)$ and matrix  $\bOm_m$ defined in \fr{eq:bOm}
satisfy the following conditions

\begin{itemize}

\item[(A1)] There exists an integer $r \geq 1$ such that
\be \label{k_cond}
\left. \frac{d^j g(t)}{dt^j} \right|_{t=0} =   \lfi
\begin{array}{ll}
0, & \mbox{if}\ \ j=0, ..., r-2,\\
B_r \ne 0, &  \mbox{if}\ \ j=r-1.
\end{array} \right.
\ee

\item[(A2)]  $g \in L_1 [0, \infty)$ is $r$ times differentiable with   $g^{(r)} \in L_1 [0, \infty)$.

\item[(A3)] Laplace transform $G(s)$ of $g$ has no zeros with nonnegative real parts except for zeros of the form $s=\infty + ib$.

\item[(A4)]  There exists $n_0$   such that for $n>n_0$, eigenvalues of matrix  $\bOm_m$  are uniformly
bounded, i.e.
\be  \label{eigenvalues}
 0< \lambda_1 \leq \lambda_{\min}(\bOm_m) \leq \lambda_{\max}(\bOm_m) \leq \lambda_2 < \infty
\ee
for some absolute constants $\lambda_1$ and $\lambda_2$.

\end{itemize}

\subsection{Introduction to theory of banded Toeplitz matrices}

The proof of asymptotic optimality of the estimator $\hat{f}_{\widehat{m}}$
relies heavily on the theory of banded Toeplitz matrices developed in
B\"{o}ttcher  and Grudsky  (2000, 2005). In this subsection, we review some of the
facts about Toeplitz matrices which we shall use later.

Consider a sequence of numbers $\{ b_k \}_{k=-\infty}^\infty$  such that $\sum_{k=-\infty}^\infty |b_k| < \infty$. 
An infinite Toeplitz matrix $T=T(b)$ is the matrix with elements $T_{i,j} = b_{i-j}$, $i,j=0,1, \ldots $.

Let ${\cal C} = \{z \in C: |z|=1 \}$ be the complex unit circle.
With each Toeplitz matrix $T(b)$ we can associate its symbol
\be \label{assocfun}
b(z) = \sum_{k=-\infty}^\infty b_k z^k, \ \ z \in {\cal C}.
\ee
Since, $\displaystyle{B(\theta)= b(e^{i\theta}) = \sum_{k=-\infty}^\infty b_k e^{i k \theta}}$,
numbers $b_k$ are Fourier coefficients of function $B(\theta)= b(e^{i\theta})$.

There is a very strong link between properties of a Toeplitz matrix $T(b)$
and function $b(z)$. In particular, if  $b(z) \neq 0$ for $z \in {\cal C}$
and $\wind (b) = J_b$, then $b(z)$ allows Wiener-Hopf factorization
$b(z) = b_{-} (z)\,  b_{+} (z)\,  z^{J_b}$  where   $b_+$ and $b_-$ have the following forms 
$$
b_{-} (z)  =  \sum_{k=0}^\infty b^{-}_{-k} z^{-k}, \ \
b_{+} (z)  =  \sum_{k=0}^\infty b^{+}_{k} z^{k}
$$
(see Theorem 1.8 of B\"{o}ttcher  and Grudsky  (2005)).

If $T(b)$ is a lower triangular Toeplitz matrix, then
$b(z) \equiv b_{+} (z)$ with $b^+_k=b_k$. 
In this case,  the product of two Toeplitz matrices can be obtained by simply multiplying their symbols and
the inverse of a Toeplitz matrix can be obtained  by taking the  reciprocal
of function  $b_{+} (z)$:
\be \label{identity}
T(b_{+} d_{+}) = T(b_{+}) T(d_{+}),\ \ \
T^{-1}(b_{+}) = T(1/b_{+}).
\ee

Let $T_m (b) = T_m (b_{+}) \in R^{m \times m}$ be a banded lower triangular
Toeplitz matrix corresponding to the Laurent polynomial
$\displaystyle{b  (z) = \sum_{k=0}^{m-1} b_k z^k}$.

In practice, one usually use only finite, banded, Toeplitz matrices with
elements $T_{i,j}$, $i,j=0,1, \ldots, m-1$.  In this case,  only a finite  number
of coefficients $b_k$ do not vanish and function $b(z)$ in \fr{assocfun}
reduces to a  Laurent polynomial $\displaystyle{b(z) = \sum_{k=-J}^{K} b_k z^k}$,
 $z \in {\cal C}$,  where $J$ and $K$ are nonnegative integers, $b_{-J} \neq 0$ and  $b_{K} \neq 0$.
If $b(z) \neq 0$ for $z \in {\cal C}$, then $b(z)$ can be represented in a form
\be \label{Laurent_pol}
b(z) = z^{-J} b_{K} \prod_{j=1}^{J_0} (z-\mu_j) \prod_{k=1}^{K_0} (z - \nu_k) \ \ \mbox{with}\ \
|\mu_j|<1,\, |\nu_k|>1.
\ee
In this case, the winding number of $b(z)$ is $\wind (b) = J_0 - J$.

Let $T_m (b) = T_m (b_{+}) \in R^{m \times m}$ be a banded lower triangular
Toeplitz matrix corresponding to the Laurent polynomial
$\displaystyle{b  (z) = \sum_{k=0}^{m-1} b_k z^k}$.
If $b$ has no zeros on the complex unit circle ${\cal C}$   and $\wind(b) =0$, then, due to Theorem 3.7
 of  B\"{o}ttcher  and Grudsky  (2005),  $T(b)$ is invertible and
$\displaystyle{\lim_{m \rightarrow \infty} \sup \|T_m^{-1} (b)\| < \infty}$.
Moreover, by Corollary 3.8,
\be \label{norm_converg}
\lim_{m \rightarrow \infty} \|T_m^{-1} (b)\| = \|T^{-1} (b)\|
\ee

\subsection{Relation between $\rho_m^2$ and $v_m^2$ }

In order to apply the theory surveyed above, we first need to
examine function $b(z)$ associated with the infinite lower triangular Toeplitz matrix
$\bG$  defined by \fr{inf_Toeplitz} and
the Laurent polynomial associated with its banded version
$\bG_m$. It turns out that $b(z)$ can be expressed via the Laplace transform $G(s)$
of the kernel $g(t)$. In particular, the following statement holds.

\begin{lemma} \label{lem:Toeplitz1}
Consider a sequence $\{ b_k \}_{k=0}^\infty$ with elements $b_0=g^{(0)}$
and $b_k = g^{(k)}  - g^{(k-1)}$, $k=1,2, \ldots$ where $g^{(k)}$ are Laguerre  coefficients of
the kernel $g$ in \fr{eq:model}. Then,   $b_k$, $k \geq 0$,  are Fourier coefficients of the function
\be \label{eq:Fourier_tr}
b(e^{i \theta}) =   G \lkr \frac{a(1+e^{i \theta})}{(1-e^{i \theta})} \rkr = \sum_{k=0}^\infty b_k e^{i \theta k},
\ee
where $G(s)$ is the Laplace transform of the kernel $g(x)$.
\end{lemma}

For any function $w(z)$ with an argument on a unit circle ${\cal C}$ denote
$$
\|w\|_{circ} = \displaystyle{\max_{|z|=1} w(z)}.
$$
The following lemma shows that indeed  $\rho_m^2  \log m  = o(v_m^2)$   as
$m \rightarrow \infty$.

\begin{lemma} \label{lem:rho_v}
Let $b(z)$ be given by \fr{eq:Fourier_tr}, i.e., $b(z) = G (a(1+z)/(1-z))$, $\|z\|=1.$
Denote
\be \label{eq:q}
w(z) = (1-z)^{-r} b(z),\ \ w^{-1} (z) = (1-z)^{r} b^{-1}(z),\ \
\|z\|=1.
\ee
Then, under assumptions (A1)--(A4), $w(z)$ and $w^{-1} (z)$ have  no zero on the complex
unit circle and,  for $m$   large enough,  one has
\begin{eqnarray}
\frac{C_r}{2 \lambda_1} \lkr \|w\|_{circ} \rkr^{-1}\  m^{2r+1} &\leq& v_m^2 \leq 2 C_r \lambda_2\, \|w^{-1} \|_{circ}\   m^{2r+1},
\label{eq:vm_norm}\\
m\, \rho_m^2   &\leq& C(r,w)\  v_m^2,  \label{eq:rhom_vm_ratio}
\end{eqnarray}
where $\rho_m^2$ and $v_m^2$ are defined in \fr{eq:vm_rhom}, $\lambda_1$ and $\lambda_2$ are given by \fr{eigenvalues}  and
$C(r,w)$ is an absolute constant which depends only on $w$ and $r$:
$$
C(r,w) = 2^{4r+1}\, [(r-1)!]^2\, \lkr \|w\|_{circ}\, \|w^{-1} \|_{circ} \rkr^2\ \lambda_2/ \lambda_1.
$$
\end{lemma}

\subsection{Asymptotic optimality of the estimators}

Note that  Lemma \ref{lem:rho_v} implies that, in   \fr{eq:penalty},
$\rho_m^2 \log m = o(v_m^2)$   as  $m \rightarrow \infty$,   so that
the second term  in \fr{eq:penalty} is  of smaller asymptotic order  than the first   term.
Consequently, as $n \rightarrow \infty$, $T/n \rightarrow 0$,  the right-hand side of \fr{eq:risk_bound}
is of the same asymptotic order as  the oracle risk \fr{eq:oracle}, so that, combination
of Theorem \ref{th:risk_bound} and Lemma \ref{lem:rho_v} leads to  the following statement.

\begin{theorem} \label{th:asymp_risk}
Let  condition \fr{eq:rho_condition}  hold for some positive constants $\alpha$ and $C_\rho$.
Then, under assumptions (A1)--(A4), for an estimator $\hat{f}_{\widehat{m}}$ of $f$ with penalty given by equation
\fr{eq:penalty} with $B>0$,    as $n \rightarrow \infty$,
\be \label{eq:optimality}
\frac{R(\hat{f}_{\widehat{m}}) }{R_{oracle}} \leq 16(1+ B) (1 + o(1)),
\ee
provided $T/n \rightarrow 0$  as $n \rightarrow \infty$.
\end{theorem}

{\bf Proof }  Let $m_0 = \arg\min_m [\|f_m - f\|^2 + \sigma^2 T n^{-1}\, v_m^2]$.
Then,  due to bounds \fr{eq:vm_norm} on  $v_m^2$, one has $ m_0 \rightarrow \infty$ and
$\frac{m_0^{2r+1} T}{n} \rightarrow 0$ as $T/n \rightarrow 0$. Hence,
it follows from Lemma \ref{lem:rho_v} that  $\rho_m^2  \log m  = o(v_m^2)$   as
$m \rightarrow \infty$ which, in combination with Theorem \ref{th:risk_bound},
completes the proof.
\\

\begin{remark} \label{rem:T}
{\em  The theory above is valid for $T$ being finite as well as for $T=T_n \rightarrow \infty$
as long as $T_n/n \rightarrow 0$ as  $n \rightarrow \infty$. Indeed, the natural consequence of $T$
being finite is that the bias term $\|f-f_m\|^2$ might be relatively large due to mis-representation
of $f$ for $t>T$. However, since both the risk
of the estimator $R(\hat{f}_{\widehat{m}})$ and the oracle risk $R_{oracle}$ are equally affected,
Theorem \ref{th:asymp_risk} remains valid whether $T=T_n$ grows with $n$ or not.}
\end{remark}

\begin{remark} \label{rem:B}
{\em  The right hand side of formula \fr{eq:optimality} is strictly increasing in $B$,
so, the smaller $B$ is, the closer the risk to the optimal oracle risk as $n \rightarrow \infty$.
Note, however, that choosing asymptotically  small value for $B$ (e.g, $B=1/n$) can make the second term  in the  penalty
 \fr{eq:penalty} dominant, so that \fr{eq:optimality} will become invalid.
}
\end{remark}

\section{Simulation study }
\label{sec:simulation}
\setcounter{equation}{0}

In order to evaluate finite sample performance of the methodology presented above, we carried out a 
simulation study. We chose three versions of the kernel $g$, normalized to have their maximum equal to 1:
\begin{itemize}
\item $g_1(t)$ which coincides with the fit of an arterial input function (AIF) for real data obtained in  
the  REMISCAN (2012) study. The real-life  observations of an AIF corresponding to kernel $g_1$ 
coming from one patient in the REMISCAN study \cite{remiscan} and   fitted estimator of $g_1$ using 
an expansion over the system of  the Laguerre functions with   $M=18$ 
are presented in Figure \ref{remiscan}. One can see clearly two behavioral patterns : 
initial high frequency behavior caused by injection of  the contrast agent  as a bolus and 
subsequent  slow decrease with regular fluctuations due to the recirculation of the contrast agent inside the blood system.

\begin{figure}[htbp]
\begin{center}
\includegraphics[width=0.45\textwidth]{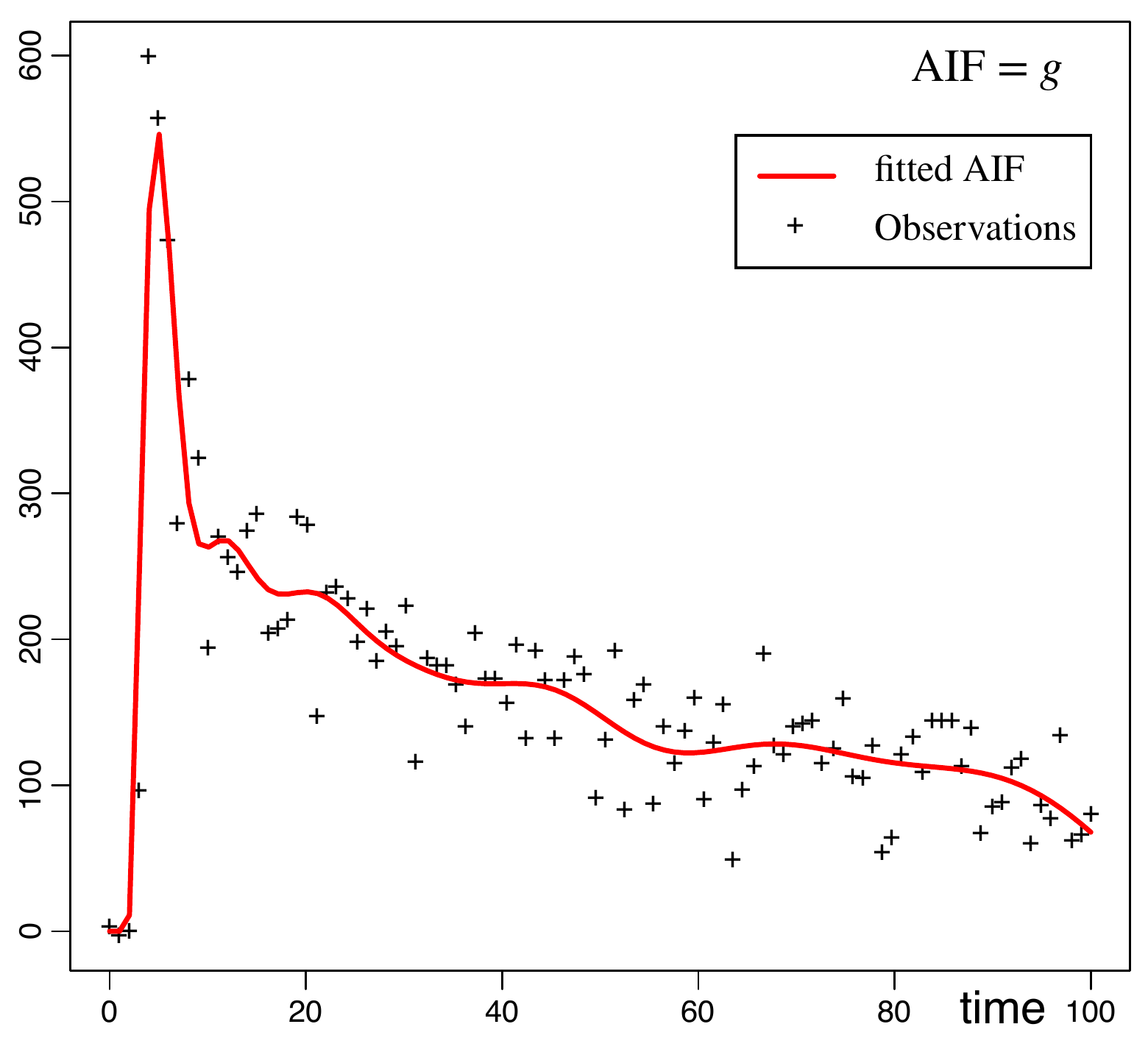}
\caption{Observations of an arterial input function (AIF) corresponding to kernel $g$ 
coming from one patient in the REMISCAN study \cite{remiscan} and fitted estimator of $g$ using 
an expansion over the system of  the Laguerre functions with   $M=17$. }
\label{remiscan}
\end{center}
\end{figure}

\item $g_2(t) =   t^2 e^{-0.1 t}$ which aims to reproduce a long injection of contrast agent; 
\item $g_3(t) = t^7 (100+t)^{-1} \, \exp \lkr -0.9 t^{3/4} \rkr$ which describes an injection with a recirculation of the contrast agent inside the blood network.
\end{itemize}
Simulations were carried out for five different test functions $f$: 
\begin{itemize}
\item $f_1(x) = \exp(-0.1 x)$, 
\item $f_2(x) = \exp(-0.6 x)$, 
\item $f_3(x) = 0.5\,  \exp(-0.1 x) + 0.5\, \exp(-0.6 x)$, 
\item $f_4(x)= 1 - IG(2; 0.5)$
where $IG(2; 0.5)$ is the cdf of the gamma distribution with the shape parameter 2 
and the scale parameter 0.5, 
\item $f_5(x) = (x+1)^{-1/3}$. 
\end{itemize}
The value of $a$ in formula \fr{eq:Laguerre_func} was chosen so that to provide the best possible fit
for the  kernel $g$ when the number of terms in the expansion of $g$ is maximum, i.e. $m = M$. 

\begin{figure}[htbp]
\begin{center}
\includegraphics[width=0.45\textwidth]{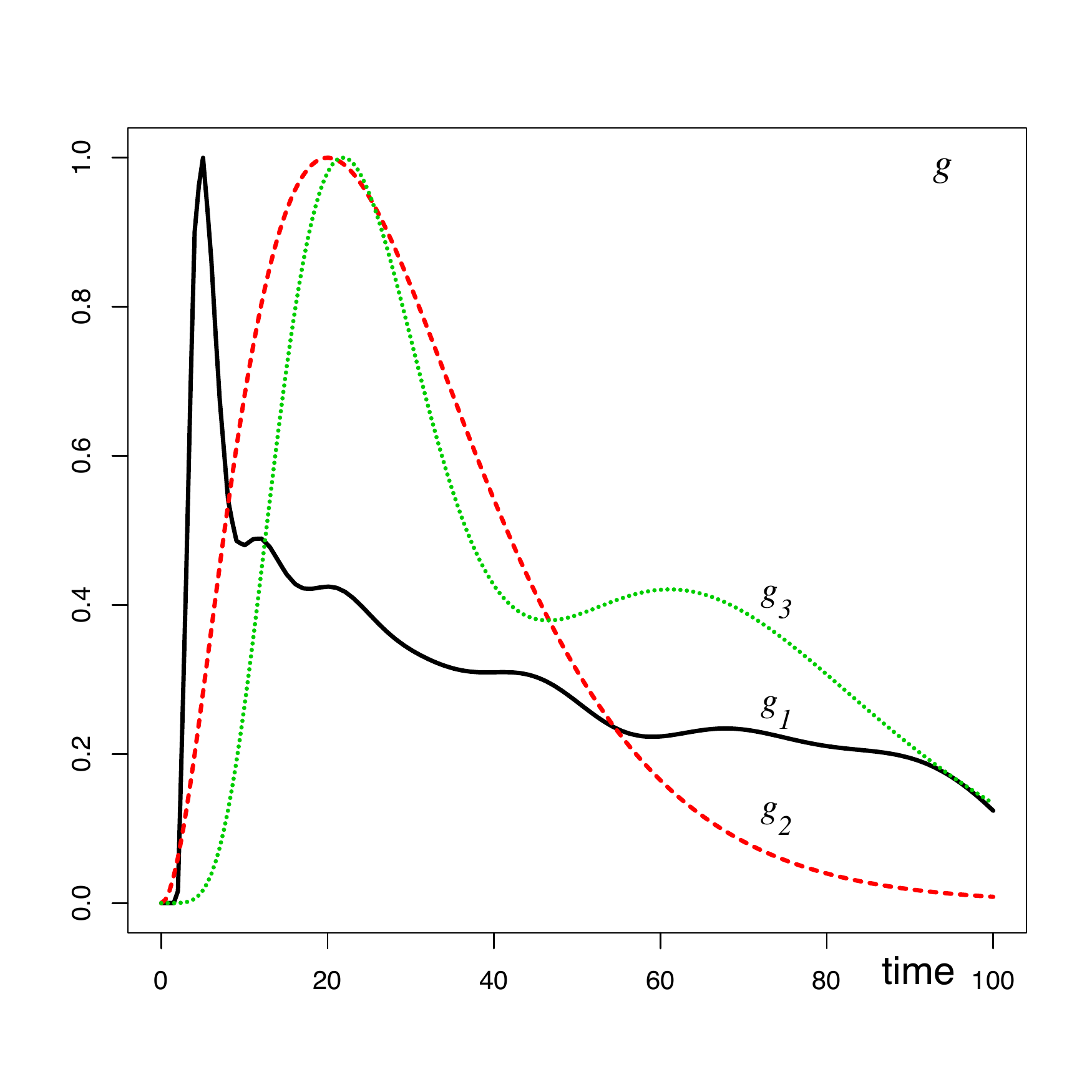}\quad \includegraphics[width=0.45\textwidth]{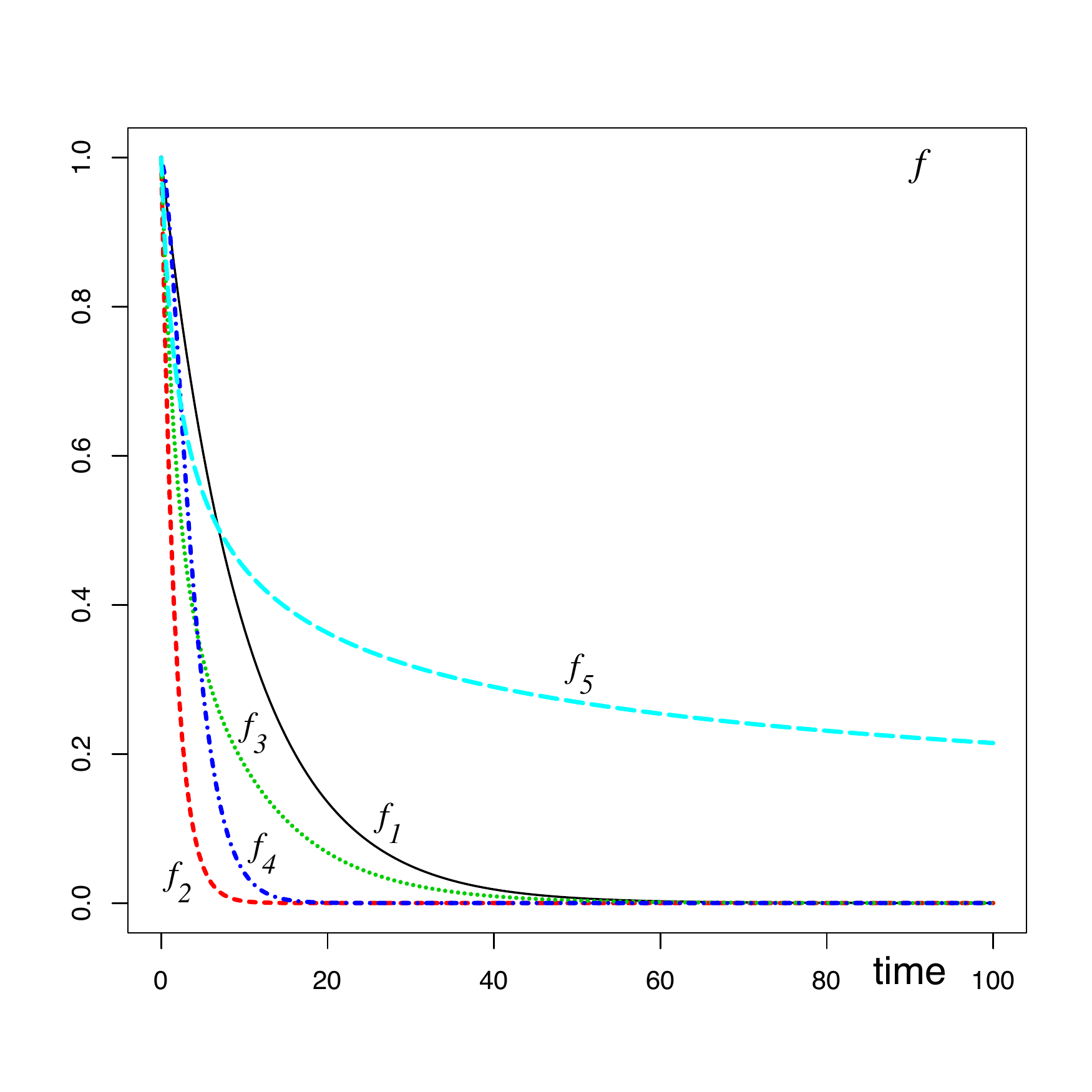}
\caption{Test functions : (left) the kernel functions $g$ - (right) the estimated functions $f$}
\label{f and g}
\end{center}
\end{figure}

The functions $f$ and $g$ are shown in Figure \ref{f and g}.

\begin{figure}[h!]
\begin{center}
\includegraphics[width=0.67\textwidth]{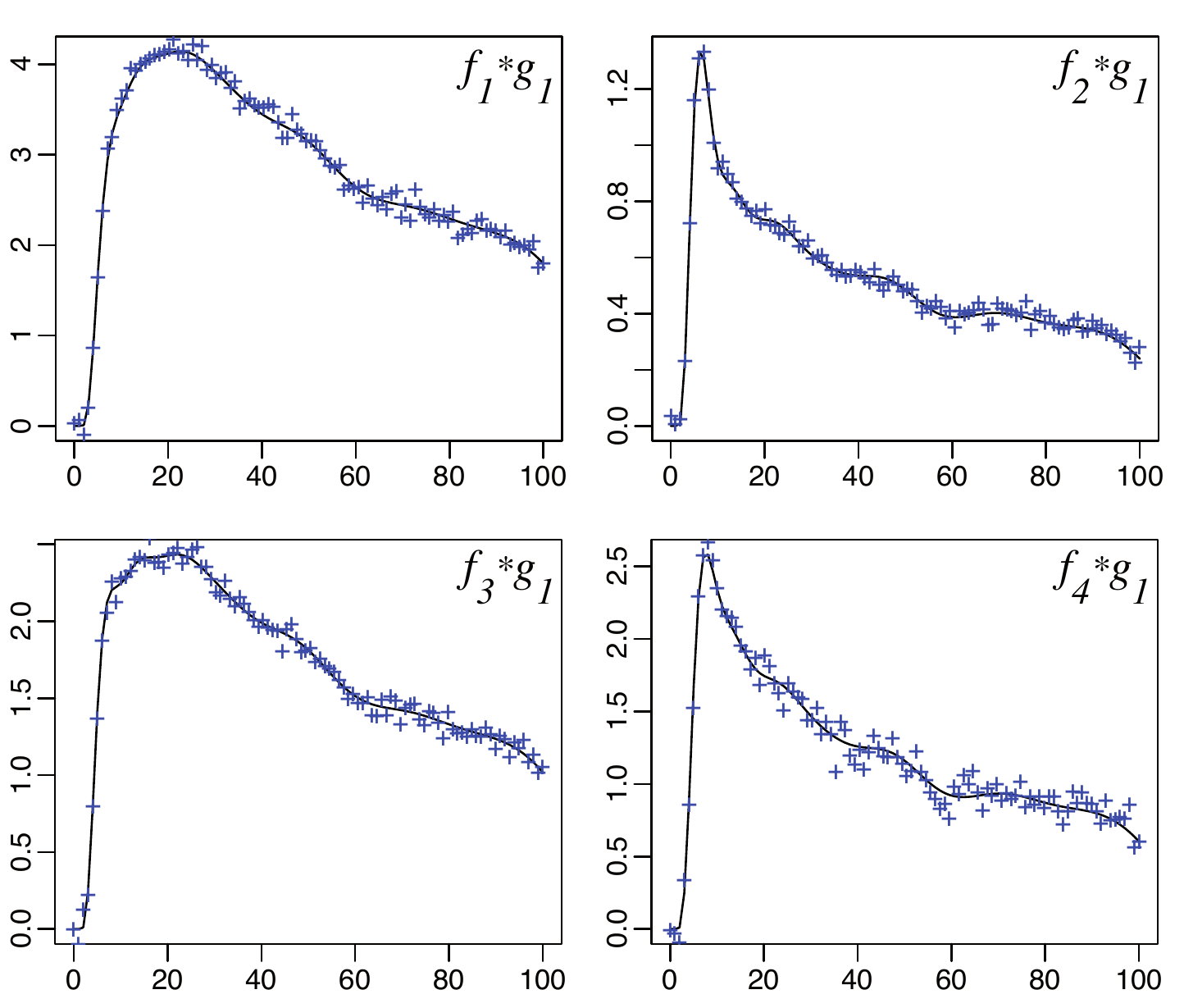}
\caption{Observations and true convolutions of kernel $g_1$ with (unknown) functions $f_1$, \ldots, $f_4$}
\label{4 examples : observations}
\end{center}
\end{figure}

\begin{figure}[!h]
\begin{center}
\includegraphics[width=0.67\textwidth]{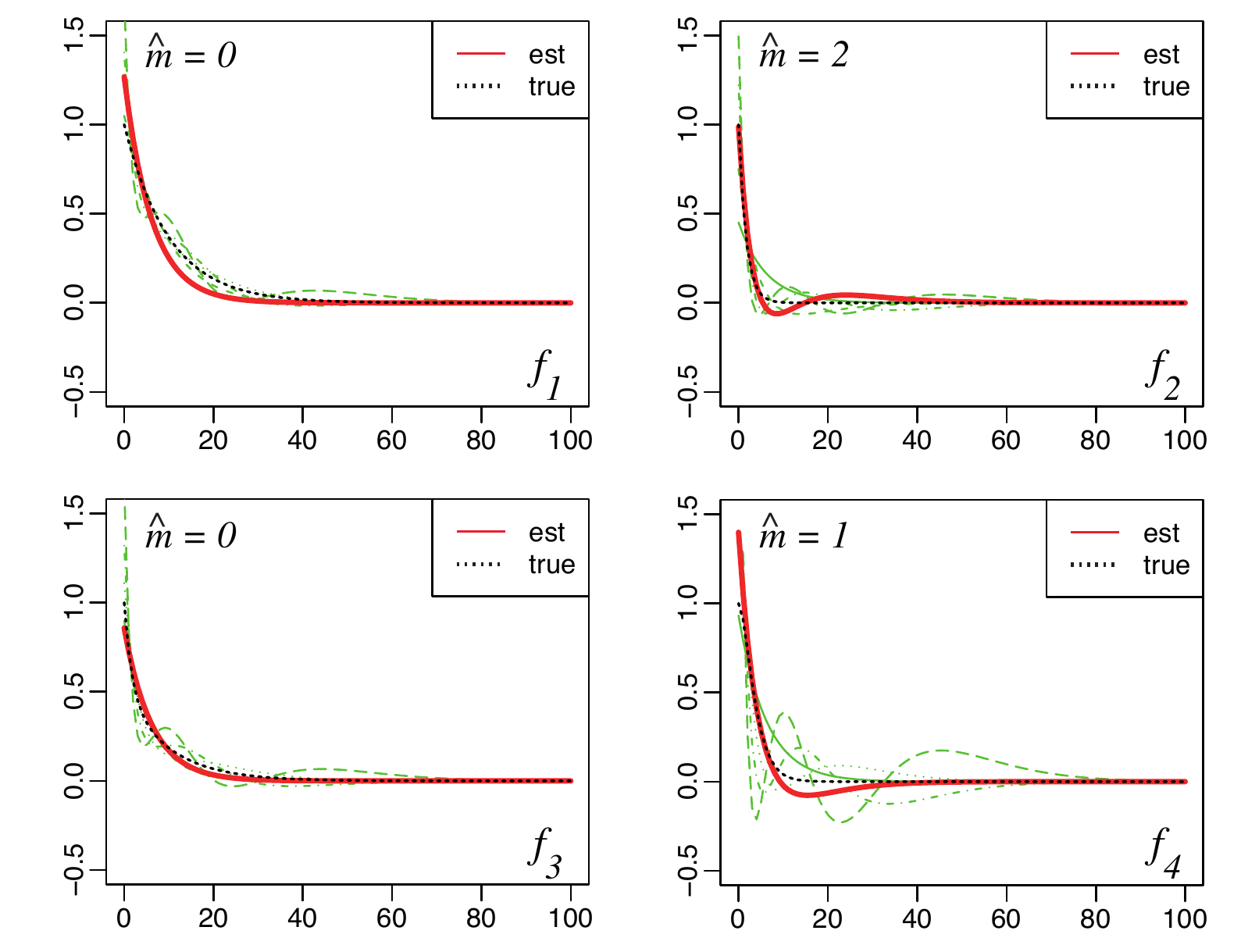}
\caption{Estimators (thick plain line) and (unknown) functions (thick dotted line) $f_1$, \ldots, $f_4$.
Other fine dashed lines represent the estimates for $m=0,\ldots,5$. The selected value of $m$ is given by $\hat{m}$.}
\label{4 examples : estimates}
\end{center}
\end{figure}

We illustrate performance of our methodology using kernel $g_1$ and test functions  $f_1$, ..., $f_4$.  
Figure \ref{4 examples : observations} shows the observations and the true convolution for a 
medium signal-to-noise ratio 8. The associated estimators are presented in Figure \ref{4 examples : estimates}. 
Here $SNR$ is defined as
$$
SNR = \sqrt{\Var(f)\big/(\sigma^2\ \Var(g))}
$$
where, for any function $\varphi$, we define $\Var (\varphi)$ as 
$$
\Var (\varphi) = \int_0^T \varphi^2(x) dx - \lkv \int_0^T \varphi (x) dx \rkv^2.
$$
The idea of defining of $SNR$ in this manner  is to remove the effect of convolution with $g$. This corresponds
to SNR of Abramovich  and Silverman   (1998).

For simulations with $g_1(t)$ we chose $\beta =1$ in \fr{eq:AIF model}.
We should mention that  the value of $\beta$ is usually unknown in real-life situations.
However, since in equation \fr{eq:AIF model}, $f(t) = \beta(1 - F(t))$ where $F(t)$ is a cdf, 
one knows that $f(0)=\beta$ and, therefore, can   estimate $\beta$  as $\hat{\beta} = \widehat{f(0)}$.

\begin{table}[!]
\caption{ The values of empirical risk increased by a factor of 100:  $100 \times \hat{R} (\hat{f})$.  
Empirical risks are computed for 400 samples for $g_1$, $g_2$ and $g_3$ and five functions of interest $f_1$, \ldots, $f_5$.}
{\scriptsize
\begin{center}
\begin{tabular}{|c|c||c|c|c|c|c||c|c|c|c|c|} \hline
\multicolumn{2}{|c||}{$100\times\hat R(\hat f)$} & \multicolumn{5}{c||}{$n=100$} &  
\multicolumn{5}{c|}{$n=200$}\\\hline
\hline
$SNR$ &    & $f_1$ & $f_2$ &  $f_3$ & $f_4$ & $f_5$ &  $f_1$ & $f_2$ &  $f_3$ & $f_4$ & $f_5$       \\\hline               
\multirow{3}{*}{5} & $g_1$ & 0.33 & 0.43 & 0.065 & 0.25 & 5.1 & 0.30 & 0.22 & 0.052 & 0.22 & 4.2\\
& $g_2$ & 0.27 & 0.37 & 0.063 & 0.25 & 5.1 & 0.22 & 0.15 & 0.051 & 0.16 & 4.0\\
& $g_3$ & 0.16 & 0.57 & 0.085 & 0.37 & 6.1 & 0.14 & 0.44 & 0.061 & 0.36 & 5.5\\\hline

\multirow{3}{*}{8} & $g_1$ & 0.31 & 0.20 & 0.061 & 0.18 & 4.0 & 0.23 & 0.12 & 0.051 & 0.070 & 3.9\\
& $g_2$ & 0.23 & 0.18 & 0.062 & 0.11 & 4.0 & 0.06 & 0.14 & 0.050 & 0.031 & 3.6\\
& $g_3$ & 0.15 & 0.57 & 0.079 & 0.37 & 5.3 & 0.13 & 0.41 & 0.059 & 0.357 & 4.8\\\hline

\multirow{3}{*}{15} & $g_1$ & 0.162 & 0.15 & 0.066 & 0.075 & 4.0 & 0.025 & 0.085 & 0.050 & 0.032 & 2.6\\
& $g_2$ & 0.022 & 0.17 & 0.061 & 0.037 & 3.4 & 0.012 & 0.115 & 0.050 & 0.027 & 2.8\\
& $g_3$ & 0.142 & 0.43 & 0.077 & 0.322 & 4.9 & 0.132 & 0.182 & 0.058 & 0.154 & 4.6\\\hline


\end{tabular}
\end{center}
}
\label{default}
\end{table}%

We executed simulations with  $T=100$, $M=11$, two values of sample sizes, $n=100$ and $n=200$, 
and three signal-to-noise ratios (SNR), namely, $SNR = 5, 8$ and 15. The value 5 corresponds to 
real-life conditions, smaller values 8 and 15 correspond to noise level attained after the first 
denoising step as described in  Rozenholc  and Rei\ss{}  (2012). 

For a given trajectory, the empirical risk  was evaluated as 
\be \label{eq:empirical_risk} 
\hat{r} (\hat{f}) = n^{-1}\ \sum_{i=1}^n \lkv \hat{f} (t_i) - f(t_i) \rkv^2
\ee
and the average empirical risk, denoted $\hat R(\hat f)$, is obtained by averaging the values of 
$\hat{r} (\hat{f})$ over  400 simulation runs.
We used  $B=1/2$ in penalty (\ref{eq:penalty})  and reduced the constant 4  in the penalty to 1.5 
since this constant is an upper bound due to a triangular inequality.  The value of $\alpha$ in 
\fr{eq:penalty} is chosen using condition \fr{eq:rho_condition}
as follows. Since $2 \log \rho_m \leq \log C_\rho + \alpha \log m$, $\alpha$ is selected by regressing   
$2 \log \rho_m$ onto $\log m$ for $m = 0, \ldots, 6$.

Results of simulations are presented in Table \ref{default}.
Table \ref{default} verifies that indeed the methodology proposed in the paper works 
exceptionally well for functions $f_i, i=1, \ldots,4$, and is still quite precise for 
test function $f_5$ for which Fourier transform   does not even  exist. 
The table demonstrates the effect of choosing parameter $a$: for function $f_3$ and $n=100$, the 
average empirical risk  does not decline when SNR grows. This is due to the bias problem arising from
the fact that $f_3$ is the sum of two exponentials and we fit only one value of $a$.
\\

\begin{figure}[htbp]
\begin{center}
\includegraphics[width=0.7\textwidth]{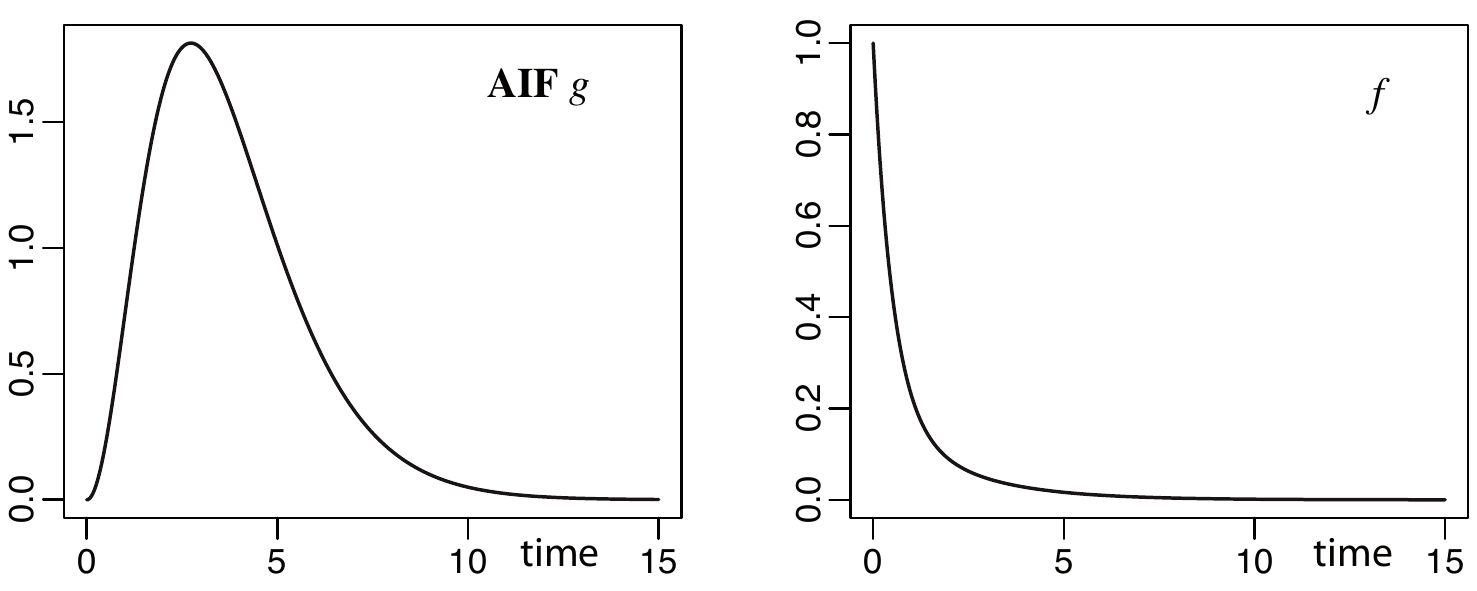}
\caption{Functions $g$ (left) and $f$ (right) for comparison between deconvolution methods using adaptive kernels and penalized Laguerre functions.}
\label{f g kernel vs laguerre}
\end{center}
\end{figure}

\begin{figure}[htbp]
\begin{center}
\includegraphics[width=0.7\textwidth]{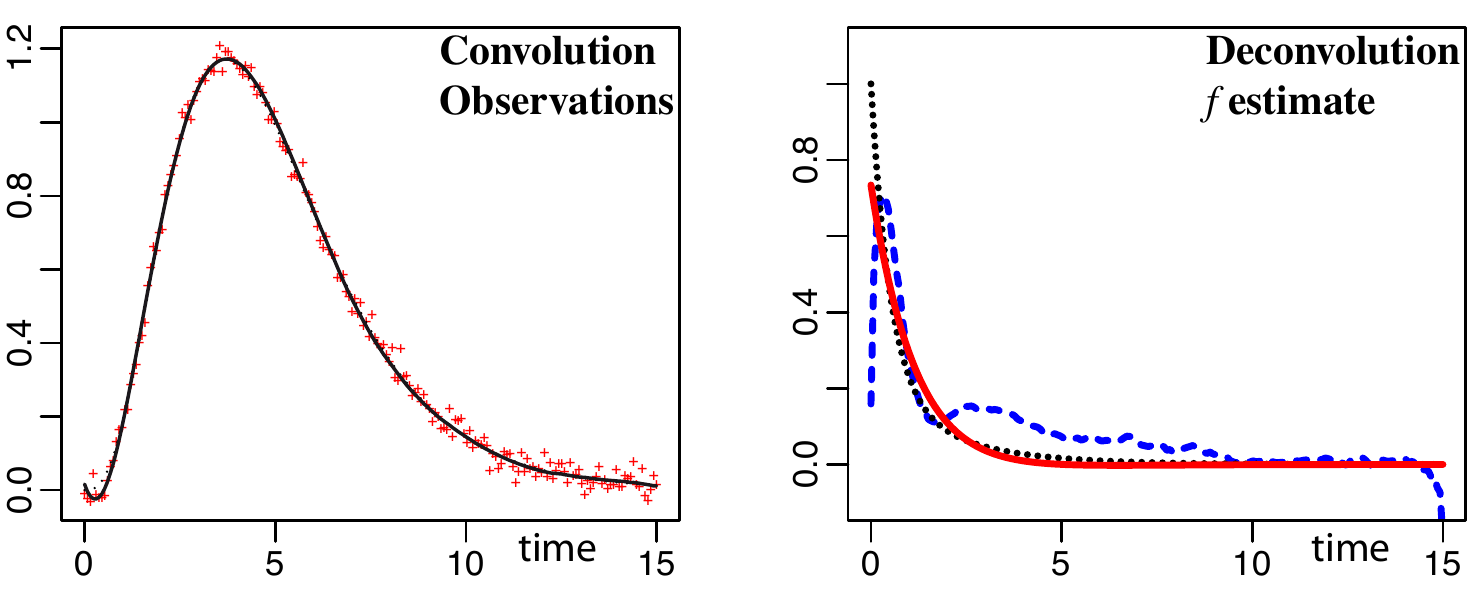}
\caption{Comparison between deconvolution using kernel method and using penalized Laguerre functions : (left) $q=f*g$ and observations; (right) estimates $f$ : Penalized Laguerre deconvolution (thick plain line) - adaptive kernel estimation (thick dashed line) - true function $f$ (dotted line).}
\label{estimate kernel vs laguerre}
\end{center}
\end{figure}

We also carried out a limited comparison of the method suggested  above with 
the technique presented in Abramovich,  Pensky and Rozenholc  (2012).
The comparison is performed using just one simple example where $f(x) = 0.2 \exp(-0.5 x) + 0.8 \exp(-2 x)$ 
and $g(t) = t^2 (t+1) e^{-t}$ (see Figure \ref{f g kernel vs laguerre}). In this example, the value of $r$ in \fr{k_cond} 
is $r=3$ and we used $n=200$, $\sigma = 0.025$ and $T=15$.
It is easy to see from Figure \ref{estimate kernel vs laguerre} that the Laguerre functions based estimator outperforms the kernel estimator of Abramovich,  Pensky and Rozenholc  (2012) and also it does not exhibit boundary effects.

\begin{figure}[htbp]
\begin{center}
\includegraphics[width=0.7\textwidth]{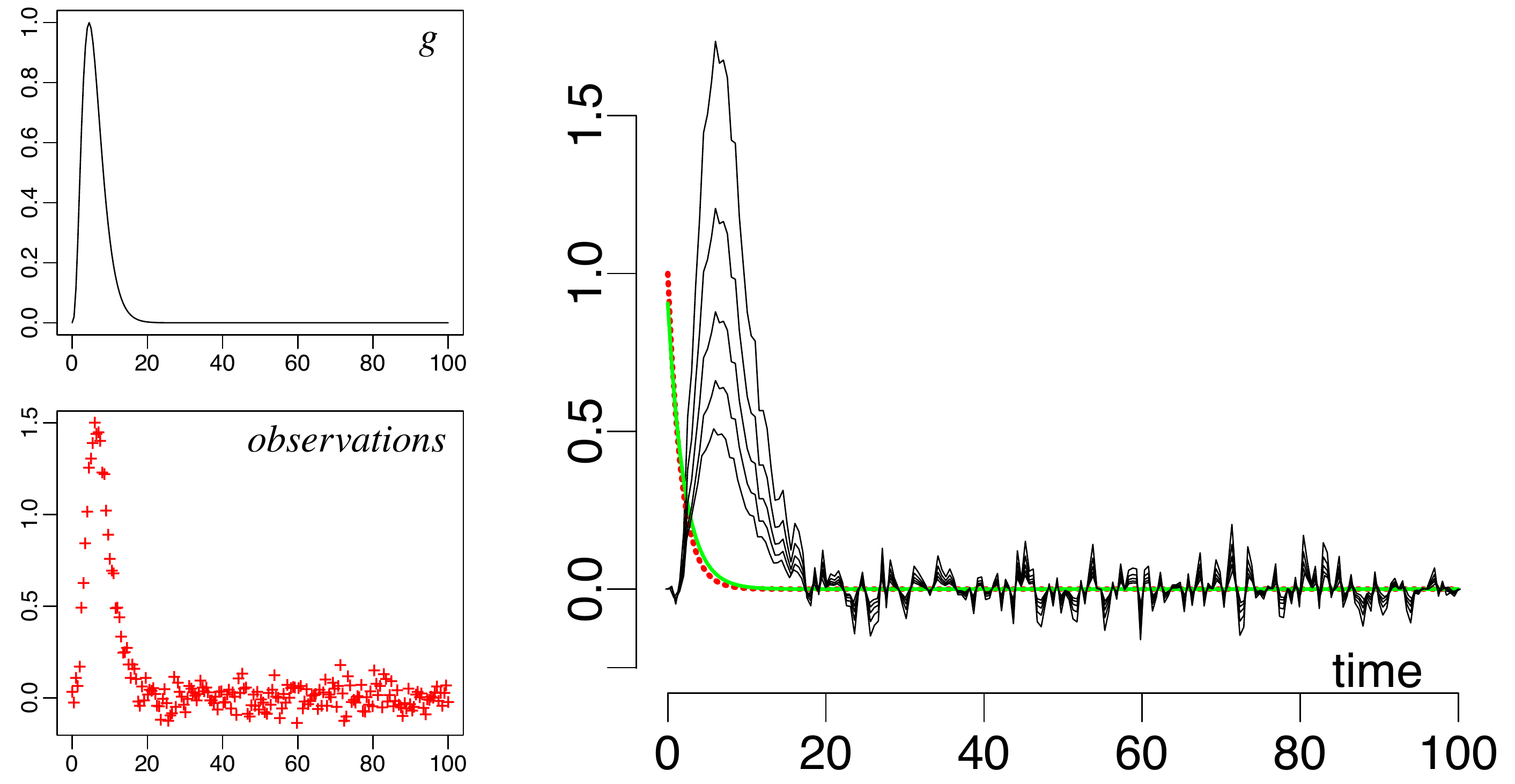}
\caption{Comparison between deconvolution using SVD method with Tikhonov regularization and penalized Laguerre functions. 
SNR=8, $n=200$, $f=f_1$, (upper-left) the kernel $g(t) = t^3 \exp(-t/3)$; (bottom-left) observations; 
(right) estimates $\hat f$ : Penalized Laguerre deconvolution (thick plain line) - SVD for various 
regularization constants (fine line) and true function $f$ (dotted line).}
\label{svd}
\end{center}
\end{figure}

Finally, we  compared our method to Singular Value Decomposition (SVD) techniques as described in the context of DCE imaging in Ostergaard {\it et al.} (1996) and Fieselmann {\it et al.} (2011). We tried various regularization
methods including thresholding and Tikhonov regularization with rectangular or trapezoid rules  for approximation 
of the convolution integral and played with the constant of regularization in order to find manually the best 
possible tuning in each case. In Figure \ref{svd} we display one of the best reconstructions which we managed to 
achieve with the  SVD approach. One can clearly see how this technique fail to adequately recover unknown function $f$ : 
first, it introduces a shift, second, it produces estimators which  fails to be a decreasing functions 
(recall that the function of interest in DCE imaging experiments is  $f(t) = \beta(1 - F(t))$  where $F(t)$ is a cdf 
and we use $\beta=1$ in our simulations).
One reason for these shortcoming is that SVD estimates are smooth and degenerate at 0. As it is noted 
 in the papers on DCE imaging (see, e.g., Fieselmann {\it et al.} (2011)), for convolution kernels corresponding to recirculation
of the contrast agent (which is a common real-life scenario), SVD fails completely and needs some extra tuning 
in order to obtain quite poor results similar to those presented in Figure \ref{svd}.

\section{Discussion }
\label{sec:discussion}
\setcounter{equation}{0}

In the present paper, we study  a noisy version of a Laplace convolution equation. 
Equations of this type frequently occur in various kinds of DCE imaging experiments. 
We propose an estimation technique for the solutions of such equation based on 
expansion of the unknown solution, the kernel and the measured right-hand side over a system 
of the Laguerre functions. The number of the terms in  the expansion  of the estimator 
is controlled via complexity penalty. The technique  leads to an estimator with the risk within a logarithmic factor 
of $m$ of the oracle risk under no assumptions on the model and within a constant factor of the oracle risk
under mild assumptions.

The major advantage of the methodology presented above is that it is usable from a practical point of view. 
Indeed, the expansion results in a small system of linear equations with the matrix of the 
system being triangular and Toeplitz. The exact knowledge of the kernel is not required: the AIF curve
can be fitted using data from DCE-CT experiments as it is shown in Figure \ref{remiscan}.
This distinguishes the present technique with the method of Abramovich,  Pensky and Rozenholc  (2012)
(referenced later as APR) which strongly depends on the knowledge of the kernel in general and the value of $r$ in \fr{k_cond}, in particular. After that, the method can be applied to any voxel of interest, either at the voxel level or using ROI
(region of interest) manually drawn by a doctor or obtained using any clustering technique.

The method is computationally very easy and fast (requires solution of a small triangular 
system of linear equations) and produces no boundary effects due to extension at zero and cut-off at $T$.
Moreover, application of the technique to discrete data does not require re-fitting the model
for each model size separately. On the contrary,  the vector of the Laguerre coefficients of the observed function
is fitted only once, for the largest model size, and then is truncated   for models of smaller sizes.
The complexity of representation of $g$ adjusts to the complexity of representation
of $f$ and the noise level. Moreover, if  $g$ can be represented by a finite expansion
over Laguerre functions with $k$ terms, the matrix of the system is $k$-diagonal.

The method performs very well in simulations. It is much more precise than the APR technique  
 as Figure \ref{estimate kernel vs laguerre} confirms. In fact, the absence of exhaustive 
comparisons between the two methods is due to the fact that it is very tricky to produce estimators 
by the APR method, especially, in the case of $g_1$ which represents real life AIF.
Similarly, as our study and Figure \ref{svd}  show, the method is much more accurate than the SVD-based
techniques.

There are few more advantages which are associated with the use of Laguerre functions basis.
Since one important goal of future analysis of DCE-CT data  is classification of the tissues 
and clustering of curves $f(t) = \beta(1-F(t))$ which characterize their blood flow properties,
representation of the curves via Laguerre basis allows to replace the problem of classification 
of curves by classification of relatively low-dimensional vectors.
In addition, due to the absence of boundary effects, the method allows to estimate 
classical medical parameters of interest  $\beta$ which describes the  perfusion of blood flow, and also
$I_f = \int   f(s)\,ds$ which characterizes the vascular mean transit time.  These parameters can be estimated 
by $\hat{\beta} = 1/\hat{f}(0)$ and $\hat{I}_f = \int \hat f(s)\,ds$, respectively.

The complexity of representation of $g$ is controlled by the choice of parameter $a$.
Parameter $a$ is a non-asymptotic constant which does not affect the convergence rates.
In practice, one can choose $a$  in order to minimize 
$\|g-\hat g_M\|$ where $\hat g_M$ is a fitted version of $g$ using the first $M$ Laguerre functions. 
Then, the same value of $a$ can be used in 
representation of the solution~$f$. Our choice of $a$ provides a reasonable trade-off
between the bias and the variance for majority of kernels considered above,  
including a real life AIF kernel coming from the REMISCAN (2007) study. 
However, our limited experimentation with choices of $a$ shows that there is room for improvement:
undeniably, fine tuning parameter $a$ can improve  estimation precision, especially, 
in the case when kernel $g$ has a strong exponential decay. However, this issue is a matter of future investigation.


\section*{Acknowledgments}

Marianna Pensky was partially supported by National Science Foundation
(NSF), grant  DMS-1106564. The authors want to express sincere gratitude to Sergei Grudski
for his invaluable help in the proof of Lemma \ref{lem:rho_v}  and very helpful discussions.


\section{Proofs }
\label{sec:proofs}
\setcounter{equation}{0}

\subsection{Proof of Theorem \ref{th:risk_bound}}
Let $m, m' \in \mathM_n$, $\bt \in S_{m'}$ and  $\bs \in S_m$.
Denote $m^* = \max(m, m') = m \vee m'$,  $\veta_m = \sigma \sqrt{T/n}\ \vxi_m$ and observe that
\begin{eqnarray}\label{obs1}
\gamma_n(\bt)-\gamma_n(\bs) &=&
\| \bt-\vf\|^2 -\|\bs-\vf\|^2 -2 \langle \bt-\bs, \bG_{M}^{-1} \veta_{m^*} \rangle,
\end{eqnarray}
where $\vf = \vf_M$ is the vector of the true $M$ first coefficients of function $f$.
Note that, due to orthonormality of the Laguerre system, for any $m$,
\begin{equation}\label{obs2}
\|\hat{f}_m  - f\|^2 = \| \widehat{\vbof}_m - \vf \|^2 + \sum_{j=M}^\infty \lkr f^{(j)}\rkr^2
\;  \mbox{ and } \; 
\|f_m  - f\|^2 = \| \vbof_m - \vf \|^2 + \sum_{j=M}^\infty \lkr f^{(j)}\rkr^2.
\end{equation}


Now, the definition of $\hat m$ yields that for any $m \in {\mathcal M}_n$ one has
$$
\gamma_n(\widehat{\vbof}_{\hat m}) + \pen (\hat m) \leq \gamma_n(\vbof_m)+  \pen (m),
$$
which with (\ref{obs1}), implies
$$\| \widehat{\vbof}_m - \vf \|^2  \leq   \| \vbof_m - \vf \|^2+ \pen(m)  + \Delta_{m,\hat m}.$$
Here
$\Delta_{m, \hat m} = 2 \langle \hat{\vf}_{\widehat{m}}- \vf_m, \bG_M^{-1} \veta_{m^\star} \rangle - \pen(\widehat{m})$, where $m^{\star} =m\vee \whm$.
Therefore, using (\ref{obs2}), we obtain that, for any $m \in {\mathcal M}_n $,
\be \label{eq:risk1}
\|\hat{f}_{\widehat{m}} - f\|^2 \leq  \|f_m - f\|^2 + \pen(m)  + \Delta_{m,\hat m}
\ee
Note that, due to $2xy\leq (x^2/4)+4 y^2$ for all $x>0, y>0$,
\beqns
 \Delta_{m,\hat m} & \leq & 2 \| \hat{\vf}_{\widehat{m}}- \vf_m \|
 \sup_{\stackrel{\bt \in S_{m^\star}} {\|\bt\|=1}} \langle \bt, \bG_M^{-1} \veta_{m^\star} \rangle - \pen(\whm) \\
& \leq & \frac 14 \| \hat{\vf}_{\widehat{m}}- \vf_m \|^2 +
4  \sup_{\stackrel{\bt \in S_{m^\star}} {\|\bt\|=1}} \langle \bt, \bG_M^{-1} \veta_{m^\star} \rangle^2
- \pen(\whm)
\eeqns

Now, denote
\be \label{eq:tau}
\tau(m, m') = \frac{\sigma^2 T}{n} \lkv (1+B) v^2_{m^*} + 2(1 + B^{-1}) (\alpha+1) \log(m^*) \rho^2_{m^*} \rkv,
\ee
where $m^* = m \vee m'$.  Since,   for any $m$,
$ \| \hat{\vf}_{\widehat{m}}- \vf_m \|^2 \leq 2 \|\hat{f}_{\widehat{m}} - f\|^2  + 2 \|f_m - f\|^2$,
then
\beqn \label{eq:Deltam}
 \Delta_{m,\hat m} & \leq & \frac 12 \, \|\hat{f}_{\widehat{m}} - f\|^2   + \frac 12 \, \|f_m - f\|^2   +
4 \lkv \sup_{\bt \in S_{m \vee \whm}}  \langle \bt, \bG_M^{-1} \veta_{m^\star} \rangle^2 - \tau(m,\whm) \rkv_+ \\
& + & 4 \tau(m,\whm) - \pen(\whm). \nonumber
\eeqn
Using the fact that  $4 \tau(m,\whm) \leq \pen(m)+ \pen(\whm)$,  combining
\fr{eq:risk1}, \fr{eq:tau} and \fr{eq:Deltam}, derive
 \beqns
\|\hat{f}_{\widehat{m}} - f\|^2 & \leq &  \|f_m - f\|^2 + \pen(m) + \frac 12 \, \|\hat{f}_{\widehat{m}} - f\|^2   + \frac 12 \, \|f_m - f\|^2\\
& + &
4 \lkv \sup_{\bt \in S_{m \vee \whm}}  \langle \bt, \bG_M^{-1} \veta_{m^\star} \rangle^2 - \tau(m,\whm) \rkv_+.
\eeqns
Finally, subtracting  $\, \|\hat{f}_{\widehat{m}} - f\|^2/2$ from both sides of the last equation and multiplying both sides
by 2, obtain
\beqn \label{eq:MISE}
\|\hat{f}_{\widehat{m}} - f\|^2 & \leq & 3 \|f_m - f\|^2  + 4 \pen(m) \\
& + &
8 \lkv \sup_{\bt \in S_{m \vee \whm}}  \langle \bt, \bG_M^{-1} \veta_{m^\star} \rangle^2 - \tau(m,\whm) \rkv_+.
\nonumber
\eeqn
Hence, validity of Theorem \ref{th:risk_bound} rests on the  following lemma which will be proved later.

\begin{lemma} \label{lem:large_deviation}
Let condition \fr{eq:rho_condition}  hold for some positive constants $\alpha$ and $C_\rho$.
Then, for any $m$ and any  $B>0$, one has
$$\EE \lkv \sup_{\bt \in S_{m \vee \whm}, \|\bt\|=1}  \langle \bt, \bG_M^{-1} \veta_{m^\star} \rangle^2 - \tau(m,\whm) \rkv_+  \leq
\frac{2 C_\rho^2 \sigma^2 T}{mn} \left( 1 + \frac{1}{B} \right).
$$
\end{lemma}

Proof of  Lemma  \ref{lem:large_deviation} is given in Section \ref{sec:supplementary}.


\subsection{Proofs of Lemmas \ref{lem:Toeplitz1} and  \ref{lem:rho_v} }

\noindent{\bf Proof  of Lemma  \ref{lem:Toeplitz1}. }
To prove this statement, we shall follow the theory of Wiener-Hopf integral equations
described in Gohberg and Feldman (1974).
Denote Fourier transform of a function $p(x)$ by
$ \hat{p} (\omega) = \displaystyle{\int_{-\infty}^\infty e^{i \omega x} p(x) dx}$
and observe that
$$
\hat{\phi}_k(\om) = (-1)^k \sqrt{2a} \frac{(a + i \om)^k}{(a- i \om)^{k+1}}.
$$
Therefore, elements of the infinite Toeplitz matrix $\bG$ in \fr{inf_Toeplitz}
are generated by the sequence $b_j$, $j \geq 0$, where
\beqn
b_j & = & (2a)^{-1/2} (g^{(j)} - g^{(j-1)}) =
\frac{1}{2\pi}\ \int_{-\infty}^\infty \hat{g} (\om)
\overline{[\hat{\phi}_j(\om) - \hat{\phi}_{j-1} (\om)]} d\om \nonumber \\
& = &  \frac{a}{\pi} \int_{-\infty}^\infty \hat{g}(\om)
\lkr \frac{i \om -a}{i \om +a} \rkr^j \frac{d\om}{a^2 + \om^2},\ \ j=0,1, \ldots.
\label{Toep_seq1}
\eeqn
 Note that $|(i \om -a)/(i \om +a)|=1$, so that we can use the following substitution in the integral \fr{Toep_seq1}:
$$
\frac{i \om -a}{i \om +a}  = e^{-i \theta} \  \Longrightarrow \ \om = \frac{a(e^{i\theta} + 1)}{i(e^{i\theta} - 1)}
= \frac{a \sin \theta}{\cos \theta -1}, \ \ 0 \leq \theta \leq 2\pi.
$$
Simple calculations show that
$$
b_j  = \frac{1}{2\pi}\ \int_{0}^{2\pi} \hat{g} \lkr  \frac{a(e^{i\theta} + 1)}{i (+e^{i\theta} - 1)} \rkr
e^{-i\theta j} d\theta,
$$
so that $b_j$, $j \in \ZZ$, are Fourier coefficients of the function
$$
B(\theta) = b(e^{i \theta})=  \hat{g} \lkr  \frac{a(e^{i\theta} + 1)}{i( e^{i\theta} - 1)} \rkr.
$$
Now, let us show that $b_j=0$ for $j<0$. Indeed, if $j=-k$, $k > 0$, then
\beqns
b_j & = &   \frac{a}{\pi} \int_{-\infty}^\infty \hat{g}(\om)
\lkr \frac{i \om + a}{i \om -a} \rkr^k \frac{d\om}{a^2 + \om^2}
= \frac{a}{\pi} \int_{-\infty}^\infty \hat{g}(\om)
\lkr \frac{i (-\om) - a}{i (-\om) + a} \rkr^k \frac{d\om}{a^2 + \om^2} \\
& = &
\frac{1}{2\pi}\ \int_{-\infty}^\infty \hat{g} (\om)
\overline{[\hat{\phi}_j(-\om) - \hat{\phi}_{j-1} (-\om)]} d\om
=  \int_{-\infty}^\infty g(x) \lkv \phi_k(-x) - \phi_{k-1}(-x) \rkv  dx = 0
\eeqns
since $g(x)=0$ if $x<0$ and $\phi_k(-x) = 0$ if $x>0$.
Hence, function $B(\theta) = b(e^{i \theta})$ has only coefficients
$b_j$, $j \geq 0$, in its Fourier series.
Now, to complete the proof, one just needs to note that $G(s) = \hat{g}(i s)$
for any $s$ such that Laplace transform $G(s)$ of $g$ exists.
\\

\noindent{\bf Proof  of Lemma  \ref{lem:rho_v}. }
Let us first find upper and  lower bounds on $\| \bG_m^{-1} \|_2^2 = \Tr(\bG_m^{-T} \bG_m^{-1})$ and
  $\| \bG_m^{-1} \|^2 = \lambda_{\max} (\bG_m^{-T} \bG_m^{-1})$.
For this purpose,  examine the function
$$
b(z)  = \hat{g} \lkr  \frac{a(z + 1)}{i (z - 1)} \rkr = G \lkr  \frac{a(z + 1)}{ 1-z } \rkr,\ \ |z|=1.
$$
Denote $y = a(z + 1)/(1-z)$, so that $z = (y-a)/(y+a)$ and  $G(y) = b((y-a)/(y+a))$.

Let us show that, under Assumptions (A1)-(A4), $b(z)$ has a zero  of order $r$ at $z=1$ and
all other zeros of $b(z)$ lie outside the unit circle.

For this purpose, assume that $y=\alpha + i \beta$ is a zero of $G$, i.e. $G(\alpha + i \beta) =0$.
Simple calculus yields
$$
\left| \frac{y-a}{y+a} \right|^2 = 1 - \frac{4 \alpha a}{(\alpha +a)^2 + \beta^2},
$$
so that $|z| =  |(y-a)/(y+a)| \leq 1$ iff $\alpha \geq 0$. But, by Assumption (A3), $G(y)$ has no zeros with
nonnegative real parts, so that $\alpha <0$ and $|z| =  |(y-a)/(y+a)| > 1$. Therefore, all zeros of $b(z)$,
which correspond to finite  zeros of $G$, lie outside the complex unit circle ${\cal C}$.

 Assumptions (A1), (A2) and properties of Laplace transform imply that
$G(s) = s^{-r} (B_r + G_r(s))$ where $G_r(s)$ is the Laplace transform of
$g^{(r)} (t)$. Hence,
$$
 \lim_{Re\   s \rightarrow \infty} s^j G(s) = \lfi
\begin{array}{ll}
0, & \mbox{if}\ \ j=0, ..., r-1,\\
B_r \ne 0, &  \mbox{if}\ \ j=r,
\end{array} \right.
$$
so that $y = \infty + i \beta$ is zero of order $r$ of $G(y)$.
Since $\displaystyle{\lim_{Re\  y \rightarrow \infty}  (y-a)/(y+a) =1}$,
$b(z)$ has zero of order $r$ at $z=1$.

Then, $b(z)$ can be written as $b(z) = (1-z)^r w(z)$ where $w(z)$ is defined by formula \fr{eq:q}
and all zeros of $w(z)$ lie outside the complex unit circle. Therefore, $w(z)$ can be written as
\be \label{w_represent}
w(z) = C_w \prod_{j=1}^N (z - \zeta_j),\ \ 0 \leq N \leq \infty,\ |\zeta_j|>1,
\ee
where $C_w$ is an absolute constant. Since $b(z)$ does not contain any negative powers of $z$
in its representation, $J_0 =0$ and $J=0$ in \fr{Laurent_pol}   and, consequently, $\wind(w) =0$.
Also, by \fr{identity} and \fr{eq:q}, one has
$T^{-1}(b) = T(b^{-1})$ where $b^{-1}(z) =   w^{-1} (z) (1-z)^{-r}$.

Now, recall that $ \| \bG_m^{-1} \|_2^2 =   \|T_m(b^{-1})\|_2^2$ and
for $\| \bG_m^{-1} \|^2 =  \|T_m(b^{-1}) \|^2$.
 Using relation between Frobenius and spectral norms $\|\bA_1 \bA_2 \|_2 \leq \|\bA_1\|_2 \|\bA_2\|$
for any matrices $\bA_1$ and $\bA_2$
(see, e.g., B\"{o}ttcher  and Grudsky  (2000), page 116), obtain
\beqn
\|T_m(b^{-1})\|_2 \leq \|T_m((1-z)^{-r})\|_2  \|T_m(w^{-1})\|, & &
 \|T_m(b^{-1})\|  \leq \|T_m((1-z)^{-r})\|   \|T_m(w^{-1})\|, \label{frob}\\
 \|T_m((1-z)^{-r})\|_2 \leq \|T_m(b^{-1})\|_2  \|T_m(w)\|, & &
\|T_m((1-z)^{-r})\|  \leq \|T_m(b^{-1})\|   \|T_m(w)\|.  \label{eigen}
\eeqn

Note that (see B\"{o}ttcher  and Grudsky  (2005), page 13)
$$
\lim_{m \rightarrow \infty}  \|T_m(w^{-1})\| =  \|w^{-1} \|_{circ},\ \
\lim_{m \rightarrow \infty}  \|T_m(w )\| = \|w\|_{circ},
$$
Also, due to representation \fr{w_represent}, both    $w$ and $w^{-1}$ are bounded,
and, therefore, $0 < \|w^{-1} \|_{circ} < \infty$ and   $0 < \|w \|_{circ} < \infty$. Denote
\be \label{def:nu_f_nu_s}
\nu_f (m) = \|T_m((1-z)^{-r})\|_2,\ \ \nu_s (m) = \|T_m((1-z)^{-r})\|.
\ee
Then, it follows from \fr{norm_converg}, \fr{frob} and \fr{eigen}   that, for $m$ large enough,
\beqn
0.5 \lkr \|w\|_{circ} \rkr^{-2}\  \nu_f^2 (m) & \leq & \|T_m(b^{-1}) \|_2^2 \leq 2 \|w^{-1} \|_{circ}^2\  \nu^2_f(m),
\label{Frobenius1}\\
0.5 \lkr \|w\|_{circ} \rkr^{-2}\  \nu_s^2 (m) & \leq & \|T_m(b^{-1}) \|^2  \leq 2 \|w^{-1} \|_{circ}^2\  \nu^2_s(m).
\label{spectral1}
\eeqn

 In order to finish the proof, we need to evaluate $\nu^2_f(m)$ and  $\nu^2_s(m)$ and also to
derive a relation between $v_m^2$, $\rho_m^2$, $\|T_m(b^{-1}) \|_2^2$ and $\|T_m(b^{-1}) \|^2$.
The first task is accomplished by the following lemma.

\begin{lemma} \label{lem:nu_f_nu_s}
Let $\nu_f(m)$ and $\nu_s (m)$ be defined in \fr{def:nu_f_nu_s}.
Then,
\beqn
 2^{-(4r-1)} [(r-1)!]^{-2} m^{2r+1} & \leq & \nu^2_f(m) \leq 0.5 m^{2r+1}, \label{Frobenius2}\\
(r!)^{-2} m^{2r} & \leq & \nu^2_s(m) \leq m^{2r}. \label{spectral2}
\eeqn
\end{lemma}

Proof of Lemma \ref{lem:nu_f_nu_s}  is given in Section \ref{sec:supplementary}.
\\

Now, to complete the proof, recall that matrix $\bOm_m$ given by \fr{eq:bOm}
is   symmetric positive definite, so that there exist an orthogonal matrix $\bU_m$
and a diagonal  matrix $\bD_m$, with eigenvalues of $\bOm$ as its diagonal elements,
such that $\bOm_m = \bU_m^T \bD_m \bU_m$ and
$\bOm_m^{-1} = \bU_m^T \bD_m^{-1} \bU_m$. Hence, by \fr{eq:Qm},  \fr{eq:vm_rhom}  and Assumption (A4)
\beqns
\|T_m(b^{-1}) \|_2^2 & = &  \Tr ( \bG_m^{-1} \bG_m^{-T})  =
\Tr ( \bG_m^{-1} \bU_m^T \bD_m \bD_m^{-1} \bU_m   \bG_m^{-T})\\
& \leq &
\| \bD_m^{-1} \|   \| \sqrt{\bD_m} \bU_m \bG_m^{-1} \|_2^2 \leq
\lambda_1^{-1} \Tr (\bG_m^{-T} \bOm_m \bG_m) =  \lambda_1^{-1} v_m^2,\\
v_m^2  & = & \| \sqrt{\bD_m}\bU_m \bG_m^{-1} \|_2^2 \leq \lambda_2 \|\bG_m^{-1} \|_2^2
= \lambda_2  \|T_m(b^{-1}) \|_2^2,  \\
\rho_m^2 & = & \| \sqrt{\bD_m}\bU_m \bG_m^{-1} \|^2 \leq \lambda_2 \|\bG_m^{-1} \|^2
= \lambda_2  \|T_m(b^{-1}) \|^2,
\eeqns
so that
\be \label{Frobenius_spectral}
\rho_m^2  \leq \lambda_2  \|T_m(b^{-1}) \|^2, \ \
\lambda_1 \|T_m(b^{-1}) \|_2^2 \leq v_m^2 \leq \lambda_2 \|T_m(b^{-1}) \|_2^2.
\ee
Combination of \fr{Frobenius1} -- \fr{Frobenius_spectral} and
Lemma \ref{lem:nu_f_nu_s} complete the proof.

\subsection {Proofs of supplementary Lemmas }
\label{sec:supplementary}

\noindent{\bf Proof of Lemma  \ref{lem:large_deviation}. }

The proof of Lemma \ref{lem:large_deviation} has two steps. The first one is the application of a $\chi^2$-type deviation inequality stated in Laurent  and Massart  (2000), and  improved by Gendre (see Lemma 3.10 of Gendre  (2009)). 
The second step consists of integrating this deviation inequality.

 The $\chi^2$-inequality is formulated as follows.  Let $\bA$   be a $p\times p$ matrix $\bA\in {\mathbb M}_p({\mathbb R})$ and $\bzeta$ be a standard Gaussian vector. Denote   $v_A^2=Tr(\bA^T \; \bA)$ and   $\rho^2(\bA)=\lambda_{{\rm max}}(\bA^T \; \bA)$. Then, for any $x>0$,
\begin{equation}\label{gendre1} 
{\mathbb P}\lkr \|\bA\bzeta\|^2\geq v_A^2 + 2\sqrt{v_A^2\rho^2(\bA) x} + \rho^2(\bA) x \rkr \leq e^{-x}.
\end{equation}

Now,  recall that for $\bt\in S_m+S_{m'}=S_{m^*}$ where   $m^*=m\vee m'$, one has
$$
\langle \bt, \bG_M^{-1} \veta_{m^*} \rangle \stackrel{\mathcal L}{=} 
\frac{\sigma^2 T}n \langle \bt_{m^*}, \sqrt{\bQ_{m^*}} \bzeta_{m^*}\rangle 
$$
where $\bt_{m*}$ is the $m^*$-dimensional vector formed by  the first $m^*$ coordinates of $\bt$ 
and $\bzeta_{m^*}$ is a standard $m^*$-dimensional Gaussian vector. Moreover,  
$$
\sup_{\bt\in S_m+S_{m'}, \|\bt\|=1}  \langle \bt, \bG_{M}^{-1}\veta_{m^*}\rangle^2
= \|\sqrt{\bQ_{m^*}} \bzeta_{m^*}\|^2. 
$$
Thus, it follows from (\ref{gendre1}) that
\begin{equation}\label{gendre3} 
  {\mathbb P} \left( \frac n{T \sigma^2}\|\sqrt{\bQ_{m^*}} \bzeta_{m^*}\|^2 \geq v_{m^*}^2 +
2\sqrt{\rho_{m^*}^2v_{m^*}^2x} +\rho_{m^*}^2 x \right)\leq e^{-x}.
\end{equation}
For any $B >0$, one has $2\sqrt{\rho_{m^*}^2v_{m^*}^2x}\leq B v_{m^*}^2 +B^{-1}  \rho_{m^*}^2 x$ so that
$$
{\mathbb P}\left( \frac n{T\sigma^2}\|\sqrt{\bQ_{m^*}} \bzeta_{m^*}\|^2 \geq (1+B) v_{m^*}^2 +
(1+B^{-1}) \rho_{m^*}^2 x\right)\leq e^{-x}.
$$
Therefore, using definition \fr{eq:tau} of $\tau(m,m')$, obtain
\begin{eqnarray*} 
&&{\mathbb E}\left(\sup_{\bt\in S_m+S_{m'}, \|\bt\|=1}  \langle \bt, \bG_M^{-1} \veta_{m^*}\rangle ^2 -\tau(m,m')\right)_+ 
= {\mathbb E}\left(  \|\sqrt{\bQ_{m^*}} \bzeta_{m^*}\|^2-\tau(m,m')\right)_+ \\
&\leq & \int_0^{+\infty} {\mathbb P}\left(\|\sqrt{\bQ_{m^*}} \bzeta_{m^*}\|^2 -\frac{\sigma^2 T}n 
\lkv (1+B) v_{m^*}^2 +2(1+B^{-1})(\alpha +1) \log(m^*)\rho_{m^*}^2 \rkv \geq  \xi\right) d\xi .
\end{eqnarray*}
Changing variables
$$
 \frac{2\sigma^2T}n (\alpha+1) (1+B^{-1}) \log(m^*)\rho_{m^*}^2 +\xi= \frac{\sigma^2T}n (1+B^{-1})\rho_{m^*}^2 x
$$ 
and application of  (\ref{gendre3}) yield
\begin{eqnarray*} 
{\mathbb E}\left(\sup_{\bt\in S_m+S_{m'}, \|\bt\|=1} \langle \bt, \bG_{M}^{-1} \veta_{m^*}\rangle ^2 -\tau(m,m')\right)_+ 
&\leq & (1+B^{-1}) \frac{\rho_{m^*}^2\sigma^2 T}n \int_{2(\alpha+1)\log(m^*)}^{+\infty} e^{-x}dx \\ 
&=& (1+B^{-1})\frac{\rho_{m^*}^2\sigma^2 T}n (m^*)^{-2(\alpha+1)} \\ & \leq & C_\rho (1+B^{-1})\frac{\sigma^2 T}n (m^*)^{-2}.
\end{eqnarray*}

Recall  that $m^\star=m\vee \hat m$ and obtain
$$
\EE \lkv \sup_{\bt \in S_{m \vee \whm}}  \langle \bt, \bG_M^{-1} \veta_{m^\star} \rangle^2 - \tau(m,\whm) \rkv_+ \leq
\sum_{m' \in  \mathM_n} \EE \left(\sup_{\bt\in S_m+S_{m'}, \|\bt\|=1}
\langle \bt, \bG_M^{-1} \veta_{m\vee m'} \rangle^2 - \tau(m,m') \right)_+ 
$$ 
and 
\begin{eqnarray*} 
&& \sum_{m'\in {\mathcal M}_n} {\mathbb E}\left(\sup_{\bt\in S_m+S_{m'}, \|\bt\|=1} 
\langle \bt, \bG_{M}^{-1} \veta_{m^*}\rangle ^2 - \tau(m,m')\right)_+\leq
C_\rho (1+B^{-1})\frac{\sigma^2 T}n \sum_{m'\in {\mathcal M}_n} (m\vee m')^{-2} \\
&\leq & C_\rho (1+B^{-1})\frac{\sigma^2 T}n \left( \sum_{m'=1}^m m^{-2}  + \sum_{m'> m} (m')^{-2}\right)\\
&\leq & C_\rho (1+B^{-1})\frac{\sigma^2 T}n \left( m^{-1}  + \int_m^{+\infty} \frac{dx}{x^2} \right)=2C_\rho (1+B^{-1})\frac{\sigma^2 T}{n m},
\end{eqnarray*}
which concludes the proof. $\Box$\\

\noindent{\bf Proof of Lemma  \ref{lem:nu_f_nu_s}. }
Note that, by formula 1.110 of Gradshtein and  Ryzhik  (1980),
$$
(1-z)^{-r} = \sum_{j=0}^\infty {r+j-1 \choose j} z^j,
$$
so that, by definition of Frobenius norm,
\beqns
\|T_m((1-z)^{-r})\|^2_2 & =&  m^2 + (m-1)^2 {r \choose 1}^2 + (m-2)^2 {r+1  \choose 2}^2 +
\ldots + {r+m-2 \choose m-1}^2 \\
& = & \sum_{j=0}^{m-1} {r+j-1 \choose j}^2 (m-j)^2,\\
\|T_m((1-z)^{-r})\|^2  & =& \max_{|z|=1} \left| \sum_{j=0}^{m-1} {r+j-1 \choose j} z^j \right|
= \sum_{j=0}^{m-1} {r+j-1 \choose r-1}.
\eeqns
If $r=1$, then
$$
\sum_{j=0}^{m-1} {r+j-1 \choose j}^2 (m-j)^2 = \sum_{j=0}^{m-1}   (m-j)^2 = \frac{m(m+1)(2m+1)}{6}.
$$
If $r \geq 2$, then
$$
\frac{j^{r-1}}{(r-1)!} \leq {r+j-1 \choose j} = \frac{(r-1+1) \ldots (r-1+j)}{(r-1)!}
\leq (j+1)^{r-1},
$$
so that, for $m \geq 4$,
\beqns
 \nu_f^2 (m) & = & \|T_m((1-z)^{-r})\|^2_2   \leq   0.5\, m^{2r+1}, \\
 \nu_f^2 (m)  & \geq & \sum_{j=m/4}^{3m/4} \frac{j^{2r-2}}{[(r-1)!]^2} (m-j)^2 \geq
\frac{m^{2r+1} 2^{-(4r-1)}}{[(r-1)!]^2},
\eeqns
which proves validity of \fr{Frobenius2}. To show that \fr{spectral2} holds,
observe that, by formula 0.151.1 of Gradshtein and  Ryzhik  (1980),
$$
\sum_{j=0}^{m-1} {r+j-1 \choose r-1} = {r+m-1 \choose r},\ \    \frac{m^r}{r!} < {r+m-1 \choose r} < m^r.
$$


\noindent
Fabienne Comte\\
Sorbonne Paris Cit\'{e}\\
Universit\'{e}  Paris  Descartes,\\
 MAP5, UMR CNRS 8145, France\\
{\em fabienne.comte@parisdescartes.fr}\\

\noindent
Charles-Andr\'{e} Cuenod\\
Sorbonne Paris Cit\'{e}\\
Universit\'{e}  Paris  Descartes, PARCC\\
European Hospital George Pompidou (HEGP-APHP)\\
LRI, INSERM U970-PARCC,   France\\
{\em ca@cuenod.net}\\

\noindent
Marianna Pensky\\
Department of Mathematics \\
University of Central Florida \\
Orlando FL 32816-1353, USA \\
{\em Marianna.Pensky@ucf.edu}\\

\noindent
Yves Rozenholc\\
Sorbonne Paris Cit\'{e}\\
Universit\'{e}  Paris  Descartes,\\
 MAP5, UMR CNRS 8145, France\\
{\em yves.rozenholc@parisdescartes.fr}

\end{document}